\documentclass[english]{amsart}
\usepackage[T1]{fontenc}
\usepackage[latin9]{inputenc}
\usepackage{babel}
\usepackage{amsthm}
\usepackage{amstext}
\usepackage{amssymb}
\usepackage{esint}
\usepackage[unicode=true]
 {hyperref}

\makeatletter
\numberwithin{equation}{section}
\numberwithin{figure}{section}

\usepackage{amscd}

\newtheorem{theorem}{Theorem}[section]
\theoremstyle{plain}

\newtheorem{definition}[theorem]{Definition}
\newtheorem{example}[theorem]{Example}
\newtheorem{babstatement}[theorem]{$\bullet$}

\newtheorem{lemma}[theorem]{Lemma}

\makeatother

\begin{document}
\title[
Generalized functions and their Fourier transforms]
{A quick description for engineering students of distributions (Generalized 
Functions) and their Fourier transforms}

\author{Michael Cwikel}

\address{Department of Mathematics, Technion - Israel Institute of Technology,
Haifa 32000, Israel}

\email{mcwikel@{math.technion.ac.il}}

\maketitle
\begin{center}
\textbf{\small Second arXiv version (17 October 2018)}

\par\end{center}

\smallskip{}

\begin{abstract}
These brief lecture notes are intended mainly for undergraduate students 
in engineering or physics or mathematics who have met or will soon be 
meeting the Dirac delta function and some other objects related to it. 
These students might have already felt - or might in the near future 
feel - not entirely comfortable with the usual intuitive explanations 
about how to ``integrate'' or ``differentiate'' or take the ``Fourier 
transform'' of these objects.

These notes will reveal to these students that there is a precise and 
rigorous way, and this also means a more useful and reliable way, to 
define these objects and the operations performed upon them. This can be 
done without any prior knowledge of functional analysis or of Lebesgue 
integration. Readers of these notes are assumed to only have studied 
basic courses in linear algebra, and calculus of functions of one and 
two variables, and an introductory course about the Fourier transform of 
functions of one variable. 

Most of the results and proofs presented here are in the framework of
the \textit{tempered distributions} introduced by Laurent Schwartz. But there
are also some very brief mentions of other approaches to distributions or generalized
functions.
\end{abstract}

\section{Pre$-$Introduction.}

These notes are intended for curious and motivated students, mainly of 
engineering or physics, who have taken a first course in Fourier 
transforms, maybe, for example, as part of the course \textit{Fourier 
series and integral transforms }which I have taught many times at 
Technion--Israel Institute of Technology.

You can read this document without knowing anything about Lebesgue integration
or functional analysis. (Maybe you will never need to study these (fascinating
and challenging!) topics, which were inspired by Fourier series and Fourier
transforms and interact very meaningfully with them.)

I have written these notes because several of my colleagues in the
Technion's Department of Electrical Engineering asked me whether it
might be possible to provide their students with a more solid 
mathematical foundation for some notions, such as the delta ``function'' 
and related ideas, which play important roles in the courses (signal 
processing etc. etc.) that these students take after learning about 
Fourier transforms. Well, almost everything is possible, but only if you 
have enough time (extra hours of lectures) and energy and determination 
to do it.

As far as I can see, these notes use only material or notions that
you have presumably already met in your introductory course about the Fourier 
transform, or in earlier courses on differential
and integral calculus and linear algebra. But they do ask you to think
about these notions in some ``new'' ways, probably quite different
from the ways that you thought about them before.

I strongly recommend that you look at least at the opening pages of the
paper \cite{FJ}. Those pages will give you some interesting insights and 
motivations about the topic of these notes. If you know some basic facts 
about Banach spaces, then that paper can take you rather further than my 
notes here can do. It has a novel approach, quite different from most 
other presentations of distribution theory. It is largely motivated
by applications of that theory to various applications in engineering, but
not to some of its traditional applications to partial differential 
equations. In any case, keep it in mind for future reading.   

Please note that you can find a list of symbols and some reminders 
about some relevant basic mathematical definitions for integration etc.\ 
at the end of these notes, in Sections \ref{symbols} and \ref{reminder}
respectively.

Let us begin by recalling the definition and a few of the basic facts
about the Fourier transform, and fixing the notation that we will use 
here for it and for one of the sets on which it operates. Other basic 
facts will be recalled later, as we need them. 

There are of course very many textbooks (and of course also internet 
sites) which present basic facts about Fourier transforms. In future 
versions of this document I may explicitly list several of them. 
Meanwhile, in this version I have occasionally referred to the books 
\cite{PZe} and \cite{PZh} and only those books, for formulations of some 
of these facts. (This is because these notes were originally intended 
for students at the Technion and these books were written, by Allan 
Pinkus and Samy Zafrany, explicitly for the above mentioned course. In 
fact \cite{PZe} and \cite{PZh} are essentially the same book in English 
and in Hebrew respectively.)

\begin{definition}\label{GandFT}
As in \cite{PZe} and \cite{PZh}. 
we will let $G(\Bbb{R)}$ denote the set of all functions $f:\Bbb{R\rightarrow\Bbb{C}}$
which are piecewise continuous (by which we mean piecewise continuous
on each bounded subinterval of $\Bbb{R}$) and which are absolutely
integrable on $\Bbb{R}$. 
For each $f\in G(\Bbb{R)}$, the Fourier
transform of $f$ is the function $\widehat{f}$ defined 
by $$\widehat{f}(\omega)=\frac{1}{2\pi}\int_{-\infty}^{\infty}e^{-i\omega x}f(x)dx$$
for all $\omega\in\Bbb{R}$. 
\end{definition}

As you probably already well know, different books use different 
definitions of the Fourier transform, formulae like 
$A\int_{-\infty}^{\infty}e^{ib\omega x}f(x)dx$ for various choices of 
the (real) constants $A$ and $b$. But, to within a change of variables, these 
are all essentially the same thing. For some of the particularly important and even 
spectacular applications of the Fourier transform, e.g. to quantum 
physics, tomography, X$-$ray crystallography, etc. etc. we also need to define 
a version of it for functions of several variables. Here we will only consider
functions (and distributions) of functions of one variable. But this
should help prepare you for the many variable case if you need it later.
(The treatment of this topic in the above mentioned paper 
\cite{FJ} includes the case of several variables.) 

We will have to use the following well
known facts several times: 
\begin{equation}
\left\{ \begin{array}{l}
\text{For each }f\in G(\Bbb{R)}\text{, the function }\widehat{f}\text{ is bounded 
and continuous,}\\
\text{and satisfies }\lim_{\omega\rightarrow\pm\infty}\widehat{f}(\omega)=0.
\end{array}\right.\label{spz}
\end{equation}
It is obvious that $\widehat{f}$ is bounded. For proofs of the other
two properties see e.g.\ Theorem 3.1 on page 94 of \cite{PZe} or page 104 of \cite{PZh}. (The proof
in these two books use the very important and useful Lebesgue Dominated Convergence Theorem, which
also gives these results for a much larger class of functions than $G(\Bbb{R})$. If you only care about
functions in $G(\Bbb{R})$, then there are much less sophisticated methods from undergraduate courses which will
also give these results.)

\smallskip

\textit{\textbf{Acknowledgement:} I warmly thank Professor Adam Shwartz 
for reading an earlier version of these notes and offering a number of 
very helpful comments and corrections of misprints. I am also very 
grateful to Yoni Roll for interesting remarks and corrections of several 
other misprints. Of course I am responsible for any other errors and 
misprints which I may have introduced since then. }

\section{Introduction.}

You have probably met the Dirac delta ``function'' in physics courses.
This is the strange object which is usually denoted by $\delta$ or
by $\delta(x)$.
It is a ``function'' which equals $0$ at every point $x\ne0$,
but at $x=0$ its ``graph'' has an infinitely high ``spike''.
We are told that somehow, miraculously, the ``area'' under this ``spike'' which has ``height''
$\infty$ and ``width'' $0$ is equal to $1$. So we can somehow
``integrate'' $\delta(x)$. Before I write down its ``integral''
let us agree here that whenever we are not quite sure that some notation
that we want to write is really a well defined mathematical object,
we shall warn ourselves of our doubts by writing that notation between
quotation marks. So we expect to have {\large ``}$\int_{a}^{b}\delta(x)dx${\large ''}$=1$
for every $a$ and $b$ such that $-\infty\le a<0<b\le\infty$. If
we believe that, then it is perhaps reasonable to suppose that we
can define the ``product'' {\large ``}$\delta(x)f(x)${\large ''}
of {\large ``}$\delta(x)${\large ''} with a continuous function
$f(x)$ on $[a,b]$ and then the ``graph'' of this product is $0$
everywhere except for a spike of ``height'' {\large ``}$f(0)\cdot\infty${\large ''}
at $x=0$. If in the ``game'' we are playing, {\large ``}$\infty\cdot0''=1$
then the ``area'' of this spike should be $f(0)$ and so {\large ``}$\int_{a}^{b}\delta(x)f(x)dx${\large ''}$=f(0)$.
But there are all sorts of problems with this ``game''. For example,
how should we define {\large ``}$\int_{a}^{b}\delta(x)f(x)dx${\large ''}
if $f$ is not continuous at $0$? We are used to thinking that if
we change the value of a function at one point $x$, for example at $x=0$,
then this should not change the value of its integral on an interval
including that point. But here this principle is wrong. Also, if we
allow ourselves to multiply numbers by $\infty$ and claim that {\large ``}$\infty\cdot0''=1$
then at least one of the associative and commutative laws for multiplication
must be wrong. If not we can deduce that \textbf{all} numbers are
equal.] Here is the ``proof'': For every two numbers $p$ and $q$,
\begin{equation}\label{foolish}
p=p\cdot1=\text{{\large``}}p\cdot(\infty\cdot0)=\infty\cdot(0\cdot p)=\infty\cdot0=\infty\cdot(0\cdot q)=q\cdot(\infty\cdot0)\text{{\large''}}=q\cdot1=q.
\end{equation}

So the delta ``function'' is \textit{not} a function in the precise
sense of the word, and if we assume it is, or use the formula ``$\infty\cdot0=1$''
carelessly, then we can easily make fools of ourselves and reach incorrect
inclusions. Let us nevertheless try to guess some more properties
of the delta ``function''. By what we said before, its ``Fourier
transform'' should be a constant function. We would expect it to
be given by 
\begin{equation}
\text{{\large``}}\widehat{\delta}(\omega)\text{{\large''}}=\frac{1}{2\pi}\text{{\large``}}\int_{-\infty}^{\infty}\delta(x)e^{-ix\omega}dx\text{{\large''}}=\frac{1}{2\pi}\cdot e^{-i\cdot0\cdot\omega}=\frac{1}{2\pi}.\label{gad}
\end{equation}
 Then, if the ``inverse Fourier theorem'' is somehow true in this
new setting, it suggests that we have some way to calculate the ``integral''
{\large ``}$\int_{-\infty}^{\infty}e^{i\omega x}d\omega${\large ''}
even though the function $e^{i\omega x}$ is not an integrable function
of $\omega$ on the interval $(-\infty,\infty)$ for any choice of
(constant) $x$. Perhaps we get 
\[
\text{{\large``}}\int_{-\infty}^{\infty}e^{i\omega x}d\omega\text{{\large''}}=2\pi\text{{\large``}}\int_{-\infty}^{\infty}\frac{1}{2\pi}e^{i\omega x}d\omega\text{{\large''}}=2\pi\text{{\large``}}\int_{-\infty}^{\infty}\widehat{\delta}(\omega)e^{i\omega x}d\omega\text{{\large''}}=2\pi\text{{\large``}}\delta(x){\large''}.\text{ }
\]

Let us now consider the function $H(t)=\text{{\large``}}\int_{-\infty}^{t}\delta(x)dx\text{{\large''}}$.
This obviously has to be $H(t)=\left\{ \begin{array}{lll}
0 & , & t<0\\
1 & , & t>0
\end{array}\right.$. But what is $H(0)$? Except for our embarrassment about deciding
the value at $0$, this function $H$ is the ``Heaviside function''.
(It is denoted by $u$ or $u_{0}$ in section 4 of the fourth chapters of \cite{PZe} and of \cite{PZh}
dealing with Laplace transforms.) 
The above ``formula\textquotedbl{} $H(t)=``\int_{-\infty}^{t}\delta(x)dx''$
could lead us to ask these next crazy questions: 

$\bullet$ Does $H$ somehow have a ``derivative'', even at $0$,
despite its discontinuity at $0$? 

$\bullet$ If so, could that ``derivative'' somehow be {\large ``}$H^{\prime}=\delta${\large ''}. 

And here are some even wilder questions: 

$\bullet$ Can we also ``differentiate'' {\large ``}$\delta${\large ''}
itself? 

$\bullet$ Can we describe its derivative precisely?

It is remarkable and surprising that these strange questions have
positive answers, which can be expressed with complete mathematical
precision. I'm going to explain those answers to you here. And it's good 
that we have them. They can enable us to use objects like the delta 
``function'' in various fields of mathematics and its applications
with greater confidence and effectiveness, with less risk of making
fools of ourselves with absurd things like \eqref{foolish}.

It took quite some time for these answers to evolve. Starting in the 
late 1800's a number of physicists and mathematicians played with these 
sorts of ideas and wondered about questions like the last one. There 
were many intuitive calculations, and good reasons, sometimes from 
physics, to believe them, and other good reasons (like our comments 
above) to doubt them.

As time passed, ways were found to treat the delta ``function''
and other related objects in a precise way. One of the most successful and
now most widely used ways of doing this was developed by the great
French mathematician Laurent Schwartz beginning with his initial work
in 1944%
\footnote{Six years later this discovery earned him the Fields Medal (generally
considered to be the equivalent of the Nobel Prize for mathematicians).\\
As well as being a brilliant mathematician, Schwartz was very actively
and deeply committed to matters of conscience and human rights. See, for example, 
\\\href{http://www-history.mcs.st-and.ac.uk/history/Obits/Schwartz.html}
{http://www-history.mcs.st-and.ac.uk/history/Obits/Schwartz.html}
and also
\\\href{http://www-history.mcs.st-and.ac.uk/history/Biographies/Schwartz.html}
{http://www-history.mcs.st-and.ac.uk/history/Biographies/Schwartz.html}\\
I had the privilege of meeting him several times and am proud to have sometimes participated, 
though only in some very small ways, in assisting some of his efforts 
together with Henri Cartan and Michel Brou\'e, when their \textit{Comit\'e des Math\'ematiciens},
tirelessly campaigned to assist mathematicians in distress.
}. The development of this topic is a particularly good example of
the way science can and should develop. On the one hand we should
not be frightened to try to work with intuitive ideas, even if at
first they seem doubtful or even partly crazy. On the other hand there
should also be a parallel process of carefully examining and trying to verify
these ideas, and seeing if they can be expressed in new more precise
and more rigorous ways. If that process succeeds, it can have all kinds of
positive consequences.

In these few pages we can only give a quick glimpse of some part of
Laurent Schwartz' work. He introduced a family of new mathematical
objects which he called \textit{distributions}. They are sometimes
also called \textit{generalized functions}. In particular, for working
with the Fourier transform, he introduced a special collection of
distributions which he called \textit{tempered distributions}. These
form a vector space%
\footnote{All vector spaces in these notes are in fact vector spaces over the
complex field. For the relevant definition and further comments see
Section \ref{reminder}.%
} which is usually denoted by $\mathcal{S}^{\prime}$. All the functions
in $G(\Bbb{R)}$ are in $\mathcal{S}^{\prime}$ and many other functions,
constants, all polynomials and many other functions which are not
integrable on $\Bbb{R}$ are in $\mathcal{S}^{\prime}$. But many
of the elements in $\mathcal{S}^{\prime}$ are not functions. In particular,
there is an \textit{exact way} of defining the delta ``function''
as an element of $\mathcal{S}^{\prime}$.

\smallskip{}
The space $\mathcal{S}^{\prime}$ has many remarkable properties.
For now we will mention only two of them:

1. It turns out to be possible to define the Fourier transform of
\textit{every} element of $\mathcal{S}^{\prime}$. If the element
happens to be a function in the space $G(\Bbb{R)}$ which we introduced
in Definition \ref{GandFT} then the new definition
and the old definition of Fourier transform give the same thing. Otherwise
the new definition may sometimes give a function, (for example $\widehat{\delta}$
is a constant, just as we guessed above) and sometimes it gives something
which is not a function. But the Fourier transform of any element
of $\mathcal{S}^{\prime}$ is always an element of $\mathcal{S}^{\prime}$.

\smallskip{}
2. It is also possible to define the derivative of every element in
$\mathcal{S}^{\prime}$. If the element happens to be a differentiable
function in the usual sense of the word and it and its derivative
also satisfies some mild ``growth'' conditions, then the new and old
definition of derivative coincide. If not, then the derivative may
fail to be a function, but it will always be an element of $\mathcal{S}^{\prime}$.
In particular we can show that our guess above that {\large ``}$H^{\prime}=\delta${\large ''}
is correct, and we can give this equation an exact meaning. We can
define the $n$th derivative of $\delta$ (or of any other tempered
distribution) for every $n\in\Bbb{N}$.

\smallskip{}

Before we can define $\mathcal{S}^{\prime}$ and describe how we define
derivatives and Fourier transforms of its elements, we need quite
a number of preliminary observations and results.

Suppose $f:\Bbb{R\rightarrow\Bbb{R}}$ is a function which arises
in the ``real world''. It could, for example be a function of time
representing, for example, an audio signal, say a bird singing. Hopefully
we have a good microphone and other equipment which will enable us
to see a good approximation to the graph of $f$ on the screen of
an oscilloscope or of a computer. But can we exactly measure or know
the value of $f(t)$ at each instant $t$? Apparently we cannot. Instead
our equipment in fact measures (approximately) the averages of $f$
over very short periods of time, i.e.\
the quantities $\frac{1}{b-a}\int_{a}^{b}f(x)dx$ for very short intervals
$[a,b]$. If $f$ is a continuous function then 
\begin{equation}
f(t)=\lim_{h\rightarrow0}\frac{1}{2h}\int_{t-h}^{t+h}f(x)dx.\label{acf}
\end{equation}
 So we can get an approximation to the value of $f(t)$ for some fixed
$t$ by using the averages of $f$ on smaller and smaller intervals
containing $t$.

Of course our equipment will not allow us to take the intervals smaller
than some strictly positive number (which depends on the amount of
money we paid for our equipment and what was its year of manufacture).
But let us now leave the birds and the ``real world'' and go back
to thinking more mathematically. The previous remarks suggest that
instead of studying a function $f:\Bbb{R\rightarrow\Bbb{C}}$ directly,
we can try to study it indirectly via the numbers $\int_{a}^{b}f(x)dx$
for all values of $a$ and $b$. But do these numbers contain \textit{all}
the information about the function? Yes, they do, at least when $f$
is continuous. This follows immediately from (\ref{acf}). A closely
related observation is contained in the next theorem.

\begin{theorem} \label{rw}Suppose that $f:\Bbb{R\rightarrow\Bbb{C}}$
and $g:\Bbb{R\rightarrow\Bbb{C}}$ are two continuous functions which
satisfy 
\begin{equation}
\int_{a}^{b}f(x)dx=\int_{a}^{b}g(x)dx\text{ for all numbers }a\text{ and }b\text{ such that }a<b.\label{zdt}
\end{equation}
 Then $f(x)=g(x)$ for all $x\in\Bbb{R}$. \end{theorem}

\textit{Proof}. This follows immediately from (\ref{acf}). $\blacksquare$

\smallskip{}

If $f$ and $g$ are only piecewise continuous%
\footnote{The definition of piecewise continuity is recalled in Section \ref{reminder}.
Here we are in fact only assuming that $f$ is piecewise continuous
on every bounded subinterval of $\Bbb{R}$.%
} then we get a conclusion which slightly weaker than in Theorem \ref{rw}.
The condition (\ref{zdt}) implies that $f(x)=g(x)$ at every point
$x\in\Bbb{R}$ where $f$ and $g$ are both continuous. So $f$ and
$g$ are ``almost'' the same function. They can differ at most on
a ``small'' set of points $E$ which only has finitely many elements
in any bounded interval. So, for example, if $f$ and $g$ both happen
to be in $G(\Bbb{R)}$ then they both have the same Fourier transform.
If we want to be able to conclude that $f(x)=g(x)$ at \textit{all}
points, including the points where they have ``jumps'' (points where
their left and right one--sided limits are different), then we can
decide, for example, to only work with piecewise continuous functions
$f$ which satisfy the extra condition 
\begin{equation}
f(x)=\frac{1}{2}\left(f(x+)+f(x-)\right)\label{avpm}
\end{equation}
 at each point of discontinuity $x$. Of course (\ref{avpm}) is also
true at all other points. It is not hard to see that the formula (\ref{acf})
is true for \textit{all} $t\in\Bbb{R}$ for functions $f$ satisfying
(\ref{avpm}) for each $x\in\Bbb{R}$. It will be convenient to have
a name for the collection of all functions $f:\Bbb{R\rightarrow\Bbb{C}}$
which are piecewise continuous on every bounded subinterval of $\Bbb{R}$
and also satisfy (\ref{avpm}) for all $x\in\Bbb{R}$. Let us call
it $PC_{bj}$. (The letters come from ``Piecewise Continuous with
Balanced Jumps''.) $PC_{bj}$ is of course a vector space.

Let us rewrite this variant of the previous theorem in a slightly
different way. Suppose that $A$ is the space $PC_{bj}$ and $B$
is the set of all characteristic functions %
\footnote{The function $\chi_{[a,b]}:\Bbb{R\rightarrow\Bbb{R}}$ is defined
by $\chi_{[a,b]}=\left\{ \begin{array}{ll}
1 & \text{if }x\in[a,b]\\
0 & \text{if }x\notin[a,b]
\end{array}\right.$. See also Section \ref{symbols}.%
} $\chi$\smallskip{}
$_{[a,b]}$ of intervals $[a,b]$ for each constant $a$ and $b$
with $-\infty<a<b<\infty$.

\begin{equation}
\left\{ \begin{array}{l}
\text{If }f\text{ and }g\text{ are in }A\text{ and}\\
\int_{-\infty}^{\infty}f(x)\phi(x)dx=\int_{-\infty}^{\infty}g(x)\phi(x)dx\mathit{\ for\ all\ functions\ }\phi\mathit{\ in\ }B\mathit{,}\\
\text{then }f(x)=g(x)\text{ for all }x\in\Bbb{R}\text{.}
\end{array}\right.\label{total}
\end{equation}

\smallskip{}
We now want to consider some other rather different examples of pairs
of sets of functions $A$ and $B$ which have the same property (\ref{total}).
In general, if $A$ and $B$ are sets of functions on $\Bbb{R}$,
we will say that $B$\textit{\ is a separating class for }$A$ if,
for every $f\in A$ and every $\phi\in B$ the function $f(x)\phi(x)$
is in $G(\Bbb{R)}$ and condition (\ref{total}) holds. We can abbreviate
this terminology and say simply that $B$ \textit{separates }$A$.
The good reason behind this terminology is that, if $f$ and $g$
are two different functions in $A$, then there is at least one element
$\phi$ in $B$ which ``separates'' them, i.e. it satisfies $\int_{-\infty}^{\infty}f(x)\phi(x)dx\ne\int_{-\infty}^{\infty}g(x)\phi(x)dx$.
(Many mathematicians use another alternative terminology here and
say that $B$\textit{\ is a total set for} $A$.)

\smallskip{}
The functions in a separating class $B$ are sometimes called ``test
functions'' because in order to know
if two functions $f$ and $g$ in $A$ are equal it is enough to test
the behaviour of their integrals with all the functions $\phi$ in
$B$.

\begin{example} \label{g1}Let $A=G(\Bbb{R)\cap}PC_{bj}$ and let
$B$ be the set of functions $\phi:\Bbb{R\rightarrow\Bbb{C}}$ of
the form $\phi(x)=e^{icx}$ for some real constant $c$. Then $B$
is a separating class for $A$. To prove this, suppose that $f$ and
$g$ in $A$ satisfy $\int_{-\infty}^{\infty}f(x)\phi(x)dx=\int_{-\infty}^{\infty}g(x)\phi(x)dx$
for all $\phi$ in $B$, i.e.\
$\int_{-\infty}^{\infty}f(x)e^{icx}dx=\int_{-\infty}^{\infty}g(x)e^{icx}dx$
for all $c\in\Bbb{R}$. So $\widehat{f}(-c)=\widehat{g}(-c)$ for
all $c\in\Bbb{R}$, i.e.\ $\widehat{f}$ and $\widehat{g}$ are the
same function. If $f$ and $g$ have left and right derivatives at
every point we can apply the inverse Fourier theorem, i.e.\ the fact
(cf.\ Theorem 3.3 on page 109 of \cite{PZe} and on page 119 of \cite{PZh}) that 
\begin{equation}
\begin{array}{l}
\lim_{R\rightarrow\infty}\int_{-R}^{R}e^{ix\omega}\widehat{u}(x)dx=\frac{1}{2}\left(u(x+)+u(x-)\right)\\
\text{if }u\in G(\Bbb{R)}\text{ and }u\text{ has right and left derivatives at }x.
\end{array}\label{ift}
\end{equation}
 to deduce that $f(x)=g(x)$ for all $x\in\Bbb{R}$. (Note that here
we have used the condition (\ref{avpm}).) If
the limit $\lim_{R\rightarrow\infty}\int_{-R}^{R}\left|f(x)-g(x)\right|^{2}dx$
is finite, then we can apply Plancherel's theorem (\cite{PZe} p. 113 or \cite{PZh} p. 124) to show that 
\[
\frac{1}{2\pi}\int_{-\infty}^{\infty}\left|f(x)-g(x)\right|^{2}dx=\int_{-\infty}^{\infty}\left|\widehat{f}(\omega)-\widehat{g}(\omega)\right|^{2}d\omega=0
\]
 and this again will give us that $f(x)=g(x)$ for all $x\in\Bbb{R}$.
If neither of these conditions hold then we need a more complicated
proof, which I will not write here, but which leads to the same conclusion.
\end{example}

\begin{example} \label{g2}Let $A$ be the set of functions $f$
in $PC_{bj}$ with the additional property that, for some positive
constant $C>0$ and some positive integer $N$, 
\begin{equation}
\left|f(x)\right|\le C\left(1+|x|\right)^{N}\text{ for all }x\in\Bbb{R}\text{.}\label{mdmfa}
\end{equation}
 (We stress that the numbers $C$ and $N$ will be different for different
functions $f$ in $A$.) Let $B$ be the set of all functions $\phi:\Bbb{R\rightarrow\Bbb{C}}$
of the form $\phi(x)=e^{-i\omega x}\cdot e^{-x^{2}}$ for some constant
$\omega\in\Bbb{R}$. Again the constant $\omega$ takes different
values for different functions $\phi$ in $B$. Let us try to compare
this example with Theorem \ref{rw}. There we were trying to get all
information about a function, or audio signal, by looking at all its
averages on all ``sharp windows''. Here we are trying to get all
information about a function, or audio signal, by looking at all the
averages of its ``frequencies'' (i.e.\ Fourier transform) on all
translations of a certain ``smooth window''.

We will now show that also in this case $B$ is a separating class
for $A$. First, since $\lim_{|x|\rightarrow\infty}e^{-x^{2}}(1+|x|)^{n+2}=0$
for every constant $n\in\Bbb{N}$, it follows that $K_{n}:=\sup_{x\in\Bbb{R}}e^{-x^{2}}(1+|x|)^{n+2}$
is finite for each $n\in\Bbb{N}$. Consequently, if $f$ is any function
in $A$, then it follows from (\ref{mdmfa}) that, for any $\phi(x)=e^{-i\omega x}e^{-x^{2}}$,
\begin{eqnarray}
\left|f(x)\phi(x)\right| & = & \left|f(x)e^{-i\omega x}e^{-x^{2}}\right|\le C\cdot\frac{e^{-x^{2}}(1+|x|)^{N+2}}{(1+|x|)^{2}}\notag\\
 & \le & C\cdot K_{N}\cdot\frac{1}{(1+|x|)^{2}}\text{ for all }x\in\Bbb{R}\text{.}\label{seb}
\end{eqnarray}
 This shows that $f\phi$ is absolutely integrable and of course it
is piecewise continuous and satisfies (\ref{avpm}) for all $x$.
Now suppose that $f$ and $g$ are any two functions in $A$ which
satisfy $\int_{-\infty}^{\infty}f(x)\phi(x)dx=\int_{-\infty}^{\infty}g(x)\phi(x)dx$
for all $\phi$ in $B$. Let us define two new functions $F(x)=f(x)e^{-x^{2}}$
and $G(x)=g(x)e^{-x^{2}}$. Then the previous condition becomes $\int_{-\infty}^{\infty}F(x)e^{-i\omega x}dx=\int_{-\infty}^{\infty}G(x)e^{-i\omega x}dx$
for all $\omega\in\Bbb{R}$. By estimates similar to (\ref{seb}),
the functions $F$ and $G$ are both in $G(\Bbb{R)}$ and so the same
arguments as in Example \ref{g1} show that $F=G$ and, consequently,
$f=g$. So, indeed, $B$ separates $A$. \end{example}

The set $A$ is ``very big'' but it does not contain all functions
in $G(\Bbb{R})\cap PC_{bj}$.%
\footnote{Consider for example the following function $f$. We set $f(x)=0$
except on each of the intervals $I_{n}=\left(n-\frac{1}{2^{n}n^{2}},n+\frac{1}{2^{n}n^{2}}\right)$
for all $n=2,3,....$ On $I_{n}$ set $f(n)=2^{n}$ and extend $f$
linearly on $\left[n-\frac{1}{2^{n}n^{2}},n\right]$ and on $\left[n,n+\frac{1}{2^{n}n^{2}}\right]$.
Since $f\left(n\pm\frac{1}{2^{n}n^{2}}\right)=0$ this means that
$f\ge0$ the region under graph of $f$ on $I_{n}$ is a triangle
of area $\frac{1}{2}\cdot\frac{2}{2^{n}n^{2}}2^{n}=\frac{1}{n^{2}}$.
It is easy to see that this function $f$ is in $G(\Bbb{R)}$ and
it is even continuous, but it is not in $A$.%
} We would like to replace $A$ by a still bigger class of functions
which contains all the functions of $G(\Bbb{R)}\cap PC_{bj}$. This
bigger class will be denoted by $SG(\Bbb{R)}$. The letters $SG$
stand for ``slow growth'' and the functions in this class are sometimes
called ``functions of slow growth''. In various books the class
of functions of slow growth is defined differently, and the notation
is also different. We have chosen a definition to suit our particular
modest purposes here%
\footnote{The reasons for choosing this definition are that, for simplicity,
we do not want to use functions which are not piecewise continuous
and we do not want to use the Lebesgue integral. The Lebesgue integral
is a very interesting and powerful generalization of the Riemann integral,
but its definition is too complicated to be discussed in these notes
or in the course ``Fourier series and integral transforms''.%
}.

\begin{definition} \label{sg}We define $SG(\Bbb{R)}$ to be the
set of all functions in $PC_{bj}$ which satisfy, for some constants
$C>0$ and $N\in\Bbb{N}$, 
\begin{equation}
\int_{-R}^{R}\left|f(x)\right|dx\le C(1+R)^{N}\text{ for all }R>0.\label{vqp}
\end{equation}
 Here again the constants $C$ and $N$ will be different for different
functions in $SG(\Bbb{R)}$. \end{definition}

Obviously $SG(\Bbb{R)}$ contains $G(\Bbb{R)\cap}PC_{bj}$.

\begin{example} \label{g3}As in Example \ref{g2}, we let $B$ be
the set of all functions of the form $\phi(x)=e^{-i\omega x}\cdot e^{-x^{2}}$
for some constant $\omega\in\Bbb{R}$.Then $B$ separates $SG(\Bbb{R)}$.
The main step for proving this is to show that, for each $f$ in $SG(\Bbb{R)}$,
the function $f(x)e^{-x^{2}}$ is absolutely integrable. This is not
difficult to do, but I ask you to believe it for the moment. It will
follow immediately from the properties \ref{phiins} and \ref{abint}
which we will soon prove. After we know that $f(x)e^{-x^{2}}$ is
integrable, the rest of the proof that $B$ separates $S(\Bbb{R)}$
uses what we know about Fourier transforms of functions in $G(\Bbb{R)}$
in exactly the same way as was done in Example \ref{g2}. \end{example}

\smallskip{}

\section{The Schwartz class $\mathcal{S}$ and some of its properties.}

\smallskip{}
Our next step will be to replace the set of functions $B$ which appeared
in the last two examples by a larger set of functions. This is a family
denoted by $\mathcal{S}$ and sometimes called \textit{the Schwartz
class}, or sometimes \textit{the class of }$C^{\infty}$\textit{\ rapidly
decreasing functions}.

\begin{definition} \label{sc}The class $\mathcal{S}$ consists of
all functions $\phi:\Bbb{R\rightarrow\Bbb{C}}$ such that the derivative
of order $n$, $\phi^{(n)}(x)$ exists for every $n\in\Bbb{N}$ and
every $x\in\Bbb{R}$ and which satisfy the condition 
\begin{equation}
\lim_{x\rightarrow\pm\infty}x^{m}\phi^{(n)}(x)=0\text{ for every pair of fixed integers }m\ge0\text{ and }n\ge0.\label{lpmi}
\end{equation}
 \end{definition}

We will see that $\mathcal{S}$ plays a central role in defining
Laurent Schwartz' tempered distributions. Let us now establish a number
of useful properties of $\mathcal{S}$. We would perhaps not expect
these properties of the very nice, very smooth, very quickly decaying
functions in $\mathcal{S}$ to have any connection with treating quite
nasty ``functions'' which are sometimes not even really functions. But we
will soon see that there is, \textit{davka}\footnote{I could not resist the temptation
to use the rather spicy Hebrew and Yiddish word ``davka'' in this place. It does not
have an exact translation into English nor apparently into most (any?) other languages. It
means something approximately like ``surprisingly'' or ``on the contrary''.}, 
a very strong connection:

\begin{babstatement} \label{phiins} Each function $\phi$ of the
form $\phi(x)=e^{-i\omega x}\cdot e^{-x^{2}}$ for some constant $\omega\in\Bbb{R}$
is in $\mathcal{S}$. \end{babstatement}

This is easy to see because for each such function $\phi$, and for
each $n\in\Bbb{N}$, $\phi^{(n)}(x)=p_{n}(x)\cdot e^{-i\omega x}\cdot e^{-x^{2}}$
for some polynomial $p_{n}$. This can easily be proved by induction
on $n$. Since $\lim_{x\rightarrow\pm\infty}x^{k}e^{-x^{2}}=0$ for
all $k\ge0$, we immediately deduce that (\ref{lpmi}) holds for all
non negative integers $m$ and $n$.

\begin{babstatement} \label{stotal} $\mathcal{S}$ separates $SG(\Bbb{R)}$.
\end{babstatement}

This follows from property \ref{phiins} and Example \ref{g3}.

\begin{babstatement} \smallskip{}
\label{equivbnd} Any infinitely differentiable function $\phi:\Bbb{R\rightarrow\Bbb{C}}$
is in $\mathcal{S}$ if and only if it satisfies 
\begin{equation}
\begin{array}{l}
\text{For each pair of fixed integers }m\ge0\text{ and }n\ge0\text{ there exists}\\
\text{a constant }C_{m,n}(\phi)\text{ such that }\left|x^{m}\phi^{(n)}(x)\right|\le C_{m,n}(\phi)\text{ for all }x\in\Bbb{R}\text{.}
\end{array}\label{gts}
\end{equation}
 \end{babstatement}

\smallskip{}
Since a continuous function on $\Bbb{R}$ which has finite limits
at $-\infty$ and $\infty$ must be bounded, it is easy to see that
(\ref{lpmi}) implies (\ref{gts}). Conversely, (\ref{gts}) implies
that $\left|x^{m}\phi^{(n)}(x)\right|=\frac{1}{|x|}\left|x^{m+1}\phi^{(n)}(x)\right|\le\frac{1}{|x|}C_{m+1,n}(\phi)$
for all $x\ne0$ and all non negative integers $m$ and $n$. This
in turn implies (\ref{lpmi}).

The constants $C_{m,n}(\phi)$ in (\ref{gts}) can be chosen to be
the numbers $\sup_{x\in\Bbb{R}}\left|x^{m}\phi^{(n)}(x)\right|$.
It is not hard to see that this supremum is also a maximum. So from
here onwards, for each $\phi\in\mathcal{S}$, we shall always use
the notation 
\begin{equation}
C_{m,n}(\phi)=\sup\left\{ \left|x^{m}\phi^{(n)}(x)\right|:x\in\Bbb{R}\right\} =\max\left\{ \left|x^{m}\phi^{(n)}(x)\right|:x\in\Bbb{R}\right\} .\label{defcmn}
\end{equation}

\begin{babstatement} \smallskip{}
\label{abint}For each $f\in SG(\Bbb{R)}$ and each $\phi\in\mathcal{S}$,
the function $f\phi$ is absolutely integrable on $\Bbb{R}$. \end{babstatement}

You might think that in the proof of property \ref{abint} which we
will give now, we are worrying too much about the constants in the
inequalities. But this is because we will need these inequalities
for another purpose later:

Choose an arbitrary function $f\in SG(\Bbb{R)}$ and let $N$ and
$C$ be the constants appearing in (\ref{vqp}). Then $|\phi(x)|\le C_{0,0}(\phi)$
and $|x^{N+2}\phi(x)|\le C_{N+2,0}(\phi)$ and so 
\[
\int_{-\infty}^{\infty}\left|f(x)\phi(x)\right|dx=\int_{-1}^{1}\left|f(x)\phi(x)\right|dx+\left(\int_{-\infty}^{-1}\left|f(x)\phi(x)\right|dx+\int_{1}^{\infty}\left|f(x)\phi(x)\right|dx\right)
\]
 
\[
\le C_{0,0}(\phi)\int_{-1}^{1}\left|f(x)\right|dx+C_{N+2,0}(\phi)\left(\int_{-\infty}^{-1}\frac{|f(x)|}{|x|^{N+2}}dx+\int_{1}^{\infty}\frac{|f(x)|}{|x|^{N+2}}dx\right).
\]
 By (\ref{vqp}) we have $\int_{-1}^{1}\left|f(x)\right|dx\le C\cdot2^{N}$.
We also have 
\begin{eqnarray*}
\int_{-\infty}^{-1}\frac{|f(x)|}{|x|^{N+2}}dx+\int_{1}^{\infty}\frac{|f(x)|}{|x|^{N+2}}dx & = & \sum_{n=1}^{\infty}\left(\int_{-n-1}^{-n}\frac{|f(x)|}{|x|^{N+2}}dx+\int_{n}^{n+1}\frac{|f(x)|}{|x|^{N+2}}dx\right)\\
 & \le & \sum_{n=1}^{\infty}\frac{1}{n^{N+2}}\left(\int_{-n-1}^{-n}|f(x)|dx+\int_{n}^{n+1}|f(x)|dx\right)\\
 & \le & \sum_{n=1}^{\infty}\frac{1}{n^{N+2}}\int_{-n-1}^{n+1}\left|f(x)\right|dx.
\end{eqnarray*}
 Using (\ref{vqp}) this last expression is dominated by 
\[
\sum_{n=1}^{\infty}\frac{1}{n^{N+2}}\cdot C(1+n+1)^{N}=C\cdot\sum_{n=1}^{\infty}\frac{(n+2)^{N}}{n^{N+2}}.
\]
 The general term in this series satisfies 
\begin{equation}
\frac{(n+2)^{N}}{n^{N+2}}=\frac{(n+2)^{N}}{n^{N}}\cdot\frac{1}{n^{2}}=\left(1+\frac{2}{n}\right)^{N}\frac{1}{n^{2}}\le3^{N}\frac{1}{n^{2}}.\label{cls}
\end{equation}

Combining all the preceding estimates, we obtain that 
\begin{equation}
\int_{-\infty}^{\infty}\left|f(x)\phi(x)\right|dx\le C\cdot\left(2^{N}C_{0,0}(\phi)+C_{N+2,0}(\phi)\cdot3^{N}\cdot\sum_{n=1}^{\infty}\frac{1}{n^{2}}\right).\label{ssgitd}
\end{equation}
 Please remember that here $N$ is a CONSTANT, and $n$ is the variable
of summation. Since the series $\sum_{n=1}^{\infty}\frac{1}{n^{2}}$
converges, the estimate (\ref{ssgitd}) shows that $f\phi$ is indeed
absolutely integrable on $(-\infty,\infty)$. (Note that by proving
this we have now also established the claim made in Example \ref{g3}.)

\begin{babstatement} \label{vecsp}\smallskip{}
$\mathcal{S}$ is a vector space. \end{babstatement}

It is obvious that whenever $\phi$ and $\psi$ are functions in $\mathcal{S}$,
then $\alpha\phi+\beta\psi$ is a function in $\mathcal{S}$ for all
complex constants $\alpha$ and $\beta$. It is very easy to verify
the remaining conditions needed to show that $\mathcal{S}$ is a vector
space. (Cf. the discussion about vector spaces in Section \ref{reminder}.)

\begin{babstatement} \label{misc}For each $\phi\in\mathcal{S}$,
each integer $n\ge0$ and each polynomial $p$, the function $p(x)\phi^{(n)}(x)$
is bounded and absolutely integrable and is also in $\mathcal{S}$.
\end{babstatement}

It follows immediately from the definitions and condition (\ref{gts}),
that $p(x)\phi^{(n)}(x)$ is in $\mathcal{S}$ and is bounded. For
exactly the same reasons $(1+x^{2})p(x)\phi^{(n)}(x)$ is, among other
things, bounded in absolute value by some constant $C$. Consequently
$\left|p(x)\phi^{(n)}(x)\right|\le\frac{C}{1+x^{2}}$ showing that
$p(x)\phi^{(n)}(x)$ is absolutely integrable.

\begin{babstatement} \label{sfour}Whenever $\phi$ is in $\mathcal{S}$,
its Fourier transform $\widehat{\phi}$ is also in $\mathcal{S}$.
\end{babstatement}

Since $\phi\in G(\Bbb{R)}$ it follows (cf.\ (\ref{spz})) that $\widehat{\phi}$
is bounded and continuous. By property \ref{misc} the function $\psi(x)=-ix\phi(x)$
is also in $G(\Bbb{R)}$ and so its Fourier transform $\widehat{\psi}$
is bounded and continuous. So we can then apply another standard result,
involving differentiation through the integral sign (with the help
of the Lebesgue dominated convergence theorem), to obtain that $\widehat{\phi}$
is differentiable and its derivative is $\widehat{\psi}$, i.e.\ $\frac{d}{d\omega}\widehat{\phi}(\omega)=\frac{1}{2\pi}\int_{-\infty}^{\infty}e^{-i\omega x}(-ix)\phi(x)dx.$
Since $\psi\in\mathcal{S}$ we can repeat the same argument to show
that $\widehat{\psi}^{\prime}=\widehat{\phi}^{\prime\prime}$ is continuous
and bounded and is the Fourier transform of yet another function in
$\mathcal{S}$. In fact we can repeat this argument as often as we
wish, and obtain that:

\begin{equation}
\begin{array}{l}
\text{For each }n\in\Bbb{N},\frac{d^{n}}{d\omega^{n}}\widehat{\phi}(\omega)\text{ exists and is}\\
\text{the Fourier transform of some function in }\mathcal{S}.
\end{array}\text{ }\label{ddez}
\end{equation}
 More precisely 
\begin{equation}
\frac{d^{n}}{d\omega^{n}}\widehat{\phi}(\omega)=\frac{1}{2\pi}\int_{-\infty}^{\infty}e^{-i\omega x}(-ix)^{n}\phi(x)dx.\label{doft}
\end{equation}

Now we will recall and use another standard result: Let $\psi$ be
a differentiable function in $G(\Bbb{R)}$ such that $\psi^{\prime}$
is also in $G(\Bbb{R)}$. Then 
\begin{equation}
\widehat{\psi^{\prime}}(\omega)=i\omega\widehat{\psi}(\omega).\label{ftd}
\end{equation}
 The first step towards obtaining (\ref{ftd}) is to observe that,
for any interval $[a,b]$, we have 
\begin{equation}
\left.\left(e^{-i\omega x}\psi(x)\right)\right|_{x=a}^{b}=\int_{a}^{b}\frac{d}{dx}\left[e^{-i\omega x}\psi(x)\right]dx=\int_{a}^{b}-i\omega e^{-i\omega x}\psi(x)dx+\int_{a}^{b}e^{-i\omega x}\psi^{\prime}(x)dx.\label{ibp}
\end{equation}
 The rest of the proof of (\ref{ftd}) is almost immediate when we
also know that ${\displaystyle \lim_{x\rightarrow\pm\infty}\psi(x)=0}$.
In fact here we will be assuming much more, i.e.\ that $\psi\in\mathcal{S}$.
Since ${\displaystyle \lim_{x\rightarrow\pm\infty}\psi(x)=0}$ we
can let $a$ tend to $-\infty$ and $b$ tend to $+\infty$ in (\ref{ibp})
and divide by $2\pi$ to obtain (\ref{ftd}). The importance of this
formula for us now is that it shows that whenever $\psi\in\mathcal{S}$,
then $\omega\widehat{\psi}(x)$ is the Fourier transform of another
function in $\mathcal{S}$. Repeating this argument $m$ times for
any integer $m\in\Bbb{N}$, we can see that $\omega^{m}\widehat{\psi}(x)$
is also the Fourier transform of a function in $\mathcal{S}$ . In
particular, since $\mathcal{S\subset}G(\Bbb{R)}$ this means that
\begin{equation}
\lim_{\omega\rightarrow\pm\infty}\omega^{m}\widehat{\psi}(\omega)=0\text{ for each integer }m\ge0.\label{psinice}
\end{equation}
 For later use let us also note more explicitly that in this way we
obtain (cf.\ (\ref{ftd})) that 
\begin{equation}
\widehat{\psi^{(m)}}(\omega)=(i\omega)^{m}\widehat{\psi}(\omega)\text{ for each }\psi\in\mathcal{S}\text{ and }m\in\Bbb{N}.\label{ftod}
\end{equation}

Now, given any $\phi\in\mathcal{S}$ and any integer $n\ge0$, let
us choose $\psi$ to be the function in $\mathcal{S}$ which exists
in view of (\ref{ddez}) with the property that $\widehat{\psi}(\omega)=\frac{d^{n}}{d\omega^{n}}\widehat{\phi}(\omega)$.
Then, by (\ref{psinice}) we have that $\lim_{\omega\rightarrow\pm\infty}\omega^{m}\frac{d^{n}}{d\omega^{n}}\widehat{\phi}(\omega)=0$.
In other words, $\widehat{\phi}\in\mathcal{S}$, and we have proved
property \ref{sfour}.

\begin{babstatement} \smallskip{}
\label{compsp}For each interval $[a,b]$ there exists a function
$\phi\in\mathcal{S}$ such that $\phi(x)>0$ for all $x\in(a,b)$
and $\phi(x)=0$ for all $x\notin(a,b)$. \end{babstatement}

The construction of such functions is not completely obvious. For
details, see Subsection \ref{csp} below.

\section{The space $\mathcal{S}^{\prime}$ of tempered distributions.}

\smallskip{}
We are now ready, or almost ready, I hope, to make the big jump from
the idea of a function to the idea of a generalized function. Examples
\ref{g1},\ref{g2} and \ref{g3} and property \ref{stotal} of $\mathcal{S}$
tell us that to determine a function $f$ in the space $G(\Bbb{R)\cap}PC_{bj}$
or even in the much bigger space $SG(\Bbb{R)}$, it is enough to know
the values of the integrals $\int_{-\infty}^{\infty}f(x)\phi(x)dx$
for all or some of the ``test functions'' $\phi\in\mathcal{S}$.
So instead of thinking about the values $f(x)$ of $f$ at each point
$x\in\Bbb{R}$ we can equivalently think about the values $L_{f}(\phi)$
of a new ``function'', which we will call $L_{f}$, at each ``point''
$\phi$ in $\mathcal{S}$. These values are given by the formula 
\begin{equation}
L_{f}(\phi):=\int_{-\infty}^{\infty}f(x)\phi(x)dx\text{ for all }\phi\in\mathcal{S}.\label{deflf}
\end{equation}

The generalized functions, or tempered distributions of Laurent Schwartz
will also be ``functions'' $L$ defined on $\mathcal{S}$ instead
of on $\Bbb{R}$. I.e., for each $\phi\in\mathcal{S}$ we have a complex
number $L(\phi)$. In other words $L$ is a map from $\mathcal{S}$
to $\Bbb{C}$. In general $L(\phi)$ will \textit{not} be given by
a formula like (\ref{deflf}). But $L$ will be required to have some
special properties:

\begin{definition} \label{td}Let $L$ be a map from the space $\mathcal{S}$
to the complex numbers which satisfies the following two conditions:

(Linearity): For each $\phi$ and $\psi$ in $\mathcal{S}$ and each
$\alpha$ and $\beta$ in $\Bbb{C}$
\[
L(\alpha\phi+\beta\psi)=\alpha L(\phi)+\beta L(\psi).
\]

(A special kind of continuity): If $\left\{ \phi_{k}\right\} _{k\in\Bbb{N}}$
is a sequence of functions in $\mathcal{S}$ such that, for each pair
of constant integers $m\ge0$ and $n\ge0$, the numbers $C_{m,n}(\phi_{k})$
defined by (\ref{defcmn}) satisfy $\lim_{k\rightarrow\infty}C_{m,n}(\phi_{k})=0$,
then 
\begin{equation}
\lim_{k\rightarrow\infty}L(\phi_{k})=0.\label{pki}
\end{equation}
 Then $L$ is said to be a \textbf{\textit{tempered distribution}}.

The set of all maps $L$ with these properties is denoted by $\mathcal{S}^{\prime}$.
\end{definition}

\smallskip{}
Before we show that the delta ``function'' and other generalized
``functions'' can be obtained as elements of $\mathcal{S}^{\prime}$
let us see that $\mathcal{S}^{\prime}$ ``contains'' all the functions
of $SG(\Bbb{R)}$. More precisely, we claim that, for each function
$f\in SG(\Bbb{R)}$, the linear map $L_{f}$ defined by the formula
(\ref{deflf}) is in $\mathcal{S}^{\prime}$. The map $L_{f}$ obviously
satisfies the linearity condition of Definition \ref{td}. To show
that the continuity condition also holds we use the estimate (\ref{ssgitd})
which immediately gives that, for every $\phi\in\mathcal{S}$,

\[
|L_{f}(\phi)|\le\int_{-\infty}^{\infty}\left|f(x)\phi(x)\right|dx\le C\cdot\left(2^{N}C_{0,0}(\phi)+C_{N+2,0}(\phi)\cdot3^{N}\cdot K\right),
\]
 where $K$ is the constant $K=\sum_{n=1}^{\infty}\frac{1}{n^{2}}$,
(it equals $\pi^{2}/6$ but the exact value is not important here)
and $C$ and $N$ are also constants which only depend on our choice
of $f$. This estimate shows that if $C_{0,0}(\phi_{k})$ and $C_{N+2,0}(\phi_{k})$
both tend to $0$ as $k$ tends to $\infty$, then (\ref{pki}) holds.

The fact that $\mathcal{S}$ separates $SG(\Bbb{R)}$ shows that for
each $f\in SG(\Bbb{R)}$ there is only one map $L$ in $\mathcal{S}^{\prime}$
such that $L=L_{f}$\smallskip{}
. In fact we can also calculate the value of $f(x)$ for each $x\in\Bbb{R}$
if we know the value of $L_{f}(\phi)$ for all $\phi\in\mathcal{S}$.
We will give an exact formula for doing this in an appendix (Subsection
\ref{recover}). So we have a well defined one to one correspondence
between functions $f\in SG(\Bbb{R)}$ and the maps $L_{f}$ in $\mathcal{S}^{\prime}$.
If we identify $f$ and $L_{f}$ we can think of $SG(\Bbb{R)}$ as
a subset of $\mathcal{S}^{\prime}$.

Now let us consider the map $\delta_{0}:\mathcal{S}\rightarrow\Bbb{C}$
defined by $\delta_{0}(\phi)=\phi(0)$ for each $\phi\in\mathcal{S}$.
This map is obviously linear: $\delta_{0}(\alpha\phi+\beta\psi)=\alpha\phi(0)+\beta\psi(0)=\alpha\delta_{0}(\phi)+\beta\delta_{0}(\psi)$.
It also satisfies the second condition of Definition \ref{td} since,
obviously, $\left|\delta_{0}(\phi)\right|\le C_{0,0}(\phi)$. So $\delta_{0}$
is a tempered distribution. It is not very difficult to show that
there is \textbf{no} function $f$ in $SG(\Bbb{R)}$ such that $\delta_{0}(\phi)=\int_{-\infty}^{\infty}f(x)\phi(x)dx$.
Again we will defer the exact proof of this to Section \ref{appendix}
(See Subsection \ref{delta}). The map $\delta_{0}$ is a well defined
map which plays the role of the delta ``function'' in an exact way.
We can also define translates of $\delta_{0}$. For each constant
$a\in\Bbb{R}$, let $\delta_{a}$ be the map from $\mathcal{S}$ to
$\Bbb{C}$ defined by $\delta_{a}(\phi)=\phi(a)$. This too is a tempered
distribution and it is the well defined version of what is sometimes
denoted by ``$\delta(x-a)$''.

\smallskip{}
I said above that we can find the derivative of the delta ``function''.
That seems crazy at first, since the delta ``function'' is not even
continuous, and it is not even a function. But if the derivative also
does not have to be a function, then maybe there is hope of finding
it.

The first step towards defining derivatives of tempered distributions
is to recall the formula for integration by parts 
\begin{equation}
\int_{a}^{b}f^{\prime}(x)\phi(x)dx=\left.\left(f(x)\phi(x)\right)\right|_{x=a}^{b}-\int_{a}^{b}f(x)\phi^{\prime}(x)dx\label{pib}
\end{equation}
 which holds whenever $f$ and $\phi$ are both continuous functions
on the interval $[a,b]$ whose derivatives $f^{\prime}$ and $\phi^{\prime}$
exist at every point of $(a,b)$ and are both Riemann integrable functions
on $(a,b)$. In particular (\ref{pib}) holds for every real $a$
and $b$ if $f$ is any function in $SG(\Bbb{R)}$ which is differentiable
at every point of $\Bbb{R}$ and whose derivative $f^{\prime}$ is
also in $SG(\Bbb{R)}$ and if $\phi$ is any function in $\mathcal{S}$.

\smallskip{}
We want to choose arbitrary $f$ and $\phi$ with these properties,
and take the limit in (\ref{pib}) as $a$ tends to $-\infty$ and
$b$ tends to $+\infty$. The problem is that we cannot be sure in
advance that $\lim_{b\rightarrow\infty}$ $f(b)\phi(b)$ and $\lim_{a\rightarrow-\infty}$
$f(a)\phi(a)$ exist. If we impose some extra condition, such as requiring
$f$ to also be in the class $A$ discussed in Example \ref{g2},
then both these limits do exist and are $0$ and so we obtain that
\begin{equation}
\int_{-\infty}^{\infty}f^{\prime}(x)\phi(x)dx=-\int_{-\infty}^{\infty}f(x)\phi^{\prime}(x)dx.\label{ibpp}
\end{equation}

In fact with a little bit more care (see Subsection \ref{subseq})
we can show that (\ref{ibpp}) holds in general, for all differentiable
$f\in SG(\Bbb{R)}$ with $f^{\prime}\in SG(\Bbb{R)}$ and all $\phi\in\mathcal{S}$.

Let us now rewrite (\ref{ibpp}) using the notation of (\ref{deflf}).
It becomes 
\begin{equation}
L_{f^{\prime}}(\phi)=-L_{f}(\phi^{\prime}).\label{asdsa}
\end{equation}
 If we think of $L_{f^{\prime}}$ as being, in some sense a ``derivative''
of $L_{f}$ then this formula suggests how to define the derivative
of a general tempered distribution:

\begin{definition} \label{deriv}For each $L\in\mathcal{S}^{\prime}$
let $DL$ be the map from $\mathcal{S}$ to $\Bbb{C}$ defined by
\[
DL(\phi)=-L(\phi^{\prime})\text{ for all }\phi\in\mathcal{S}.
\]
 $DL$ is called the \textbf{\textit{derivative}} or the \textbf{\textit{distributional
derivative}} of $L$. \end{definition}

We need a few moments to check that $DL$ is also in $\mathcal{S}^{\prime}$.
Since $L$ is linear, and $(\alpha\phi+\beta\psi)^{\prime}=\alpha\phi^{\prime}+\beta\psi^{\prime}$
we immediately obtain that $DL$ is linear. The second continuity
condition is also easily verified: If $\left\{ \phi_{k}\right\} _{k\in\Bbb{N}}$
is a sequence in $\mathcal{S}$ satisfying $\lim_{k\rightarrow\infty}C_{m,n}(\phi_{k})=0$
for all non negative integers $m$ and $n$, then the sequence $\left\{ \phi_{k}^{\prime}\right\} _{k\in\Bbb{N}}$
has the same properties, simply because $C_{m,n}(\phi_{k}^{\prime})=C_{m,n+1}(\phi_{k})$.
So $\lim_{k\rightarrow\infty}DL(\phi_{k})=-\lim_{k\rightarrow\infty}L(\phi_{k}^{\prime})=0$.

We can now extend this definition to define the $n^{th}$ order derivative
of $L$ to be the map $D^{n}L$ given by

\[
D^{n}L(\phi)=(-1)^{n}L(\phi^{(n)})\text{ for all }\phi\in\mathcal{S}.
\]
 It is clear that $D^{n}L=D(D^{n-1}L)$ and that all these derivatives
are elements of $\mathcal{S}^{\prime}$.

The formula (\ref{asdsa}) tells us that whenever $f$ is a differentiable
function in $SG(\Bbb{R)}$ such that $f^{\prime}$ is also in $SG(\Bbb{R)}$
then $DL_{f}=L_{f^{\prime}}$. We also see that $D\delta_{0}$ is
given by $D\delta_{0}(\phi)=-\phi^{\prime}(0)$. More generally, for
each $a\in\Bbb{R}$ and $n\in\Bbb{N}$, the $n^{th}$ derivative $D^{n}\delta_{a}$
is defined by $D^{n}\delta_{a}(\phi)=(-1)^{n}\phi^{(n)}(a)$ for all
$\phi\in\mathcal{S}$. As our last example, consider the Heaviside
function $H$ already mentioned above and defined by $H(t)=\left\{ \begin{array}{lll}
0 & , & t<0\\
1 & , & t>0
\end{array}\right.$. It corresponds to the tempered distribution $L_{H}$ defined by
$L_{H}(\phi)=\int_{0}^{\infty}\phi(x)dx$ for all $\phi\in\mathcal{S}$.
So its derivative, $DL_{H}$ is defined by $DL_{H}(\phi)=-L_{H}(\phi^{\prime})=-\int_{0}^{\infty}\phi^{\prime}(x)dx=-(\lim_{r\rightarrow\infty}\phi(r)-\phi(0))$.
Since the limit is $0$ we obtain $DL_{H}(\phi)=\phi(0)$ for all
$\phi\in\mathcal{S}$, i.e.\ we have now shown in an exact way that
$DL_{H}=\delta_{0}$, the derivative of (the distribution corresponding
to) $H$ really is the delta ``function''.

I said in the heading to this section that $\mathcal{S}^{\prime}$
is a space. Indeed it is a vector space. Given any $L$ and $M$ in
$\mathcal{S}^{\prime}$ and any complex numbers $\alpha$ and $\beta$
we define the map $\alpha L+\beta M$ from $\mathcal{S}$ to $\Bbb{\Bbb{C}}$
in the obvious way, i.e.\ $(\alpha L+\beta M)(\phi)=\alpha L(\phi)+\beta M(\phi)$
for each $\phi\in\mathcal{S}$. It is easy to check that $\alpha L+\beta M$,
defined in this way, is also an element of $\mathcal{S}^{\prime}$.
The other conditions needed to show that $\mathcal{S}^{\prime}$ is
a vector space are easily verified%
\footnote{Cf. the general remarks about complex vector spaces 
in Section \ref{reminder}.%
}. For some remarks about the choice of the notation $\mathcal{S}^{\prime}$
for the space of tempered distributions see Subsection \ref{name}
of the appendix.

\smallskip{}
Now we shall define the Fourier transforms of tempered distributions.
The idea for doing this is in a similar spirit to what was done to
define derivatives. First we find a suitable formula which is satisfied
by Fourier transforms of functions. Then we ``translate'' this formula
into something which also makes sense when we have distributions instead
of functions.

Let us recall a formula for Fourier transforms of \textbf{\textit{functions}}
which is perhaps reminiscent of the generalized Plancherel formula (\cite{PZe}, 
p. 114 or \cite{PZh} p. 124).
Before stating it, we should emphasize that, although it is ``traditional''
to use the variable $x$ for a function $f$ and the variable $\omega$
(or sometimes $\xi$) for its Fourier transform $\widehat{f}$, we
have to be ready to change the names of these variables whenever necessary.
For example, instead of saying that the Fourier transform of $e^{-|x|}$
is $\frac{1}{\pi(\omega^{2}+1)}$ it is equivalent, and maybe more
precise to say: If $f:\Bbb{R\rightarrow\Bbb{R}}$ is the function
defined by $f(t)=e^{-|t|}$ for all $t\in\Bbb{R}$, then the Fourier
transform of $f$ is the function $\widehat{f}:\Bbb{R\rightarrow\Bbb{R}}$
defined by $\widehat{f}(t)=\frac{1}{\pi(t^{2}+1)}$ for all $t\in\Bbb{R}$.
The choice of the particular letters (here we chose $t$ for both
$f$ and $\widehat{f}$) in the process of defining the function and
its transform is completely unimportant. We could equally well choose
$x$ or $\omega$, or any other letter, in each case. In particular,
in the next result we will be writing $\widehat{f}(x)$ instead of
the more traditional $\widehat{f}(\omega)$. But the meaning should
be clear.

\begin{lemma} \label{gpf}Suppose that $f$ and $\phi$ are functions
in $G(\Bbb{R)}$. Then the product functions $f(x)\widehat{\phi}(x)$
and $\widehat{f}(x)\phi(x)$ are also in $G(\Bbb{R)}$ and 
\begin{equation}
\int_{-\infty}^{\infty}\widehat{f}(x)\phi(x)dx=\int_{-\infty}^{\infty}f(x)\widehat{\phi}(x)dx.\label{aswd}
\end{equation}
 \end{lemma}

\textit{Proof. }By (\ref{spz}) $\widehat{f}$ and $\widehat{\phi}$
are bounded and continuous. This immediately implies that $f\widehat{\phi}$
and $\widehat{f}\phi$ are both in $G(\Bbb{R})$ and so both of the
integrals in (\ref{aswd}) exist. To prove that these integrals are
equal we shall use Fubini's theorem%
\footnote{This is perhaps the only result which you have not explicitly met
before. But it is needed to justify some of the basic results for
the Fourier transform which you presumably know (the inversion theorem,
Plancherel formula and convolution property). It is used, at least
implicitly, in proofs on pages 109, 114 and 117 of \cite{PZe} 
(corresponding to pages 119, 124 and 127 of \cite{PZh}).%
} which enables us, \textbf{\textit{provided that certain conditions
are fulfilled, }}to change the order of integration in repeated integrals
on $(-\infty,\infty)$. For details about a simple special case of this theorem 
which is sufficient for our needs here see, for example, Appendix B on pages 186-7 
of \cite{PZe}, or the web document,

\href{http://www.math.technion.ac.il/~mcwikel/FUBINI.pdf}{http://www.math.technion.ac.il/~mcwikel/FUBINI.pdf}

Let us consider the function $G(x,y)=f(x)\phi(y)e^{-iyx}$ on $\Bbb{R}^{2}$:

(i) $G(x,y)$ is ``piecewise continuous'' in a certain sense, (a
product of piecewise continuous functions of one variable with a continuous
function of two variables, as specified in equation (1) of the web
document just referred to).

(ii) The limit $\lim_{N\rightarrow\infty}\int_{-N}^{N}\left(\int_{-N}^{N}|G(x,y)|dy\right)dx$
is finite. This is because

\noindent $\int_{-N}^{N}\left(\int_{-N}^{N}|G(x,y)|dy\right)dx$ increases
with $N$ and it is bounded. It is bounded because 
\[
\int_{-N}^{N}\left(\int_{-N}^{N}|G(x,y)|dy\right)dx=\int_{-N}^{N}\left(\int_{-N}^{N}|f(x)\phi(y)|dy\right)dx
\]
 
\[
=\int_{-N}^{N}|f(x)|dx\cdot\int_{-N}^{N}\left|\phi(y)\right|dy\le\int_{-\infty}^{\infty}|f(x)|dx\cdot\int_{-\infty}^{\infty}\left|\phi(y)\right|dy
\]
 and this last expression is a product of two finite quantities which
do not depend on $N$.

\smallskip{}
In view of the properties (i) and (ii) we can apply Fubini's theorem
to obtain that

\[
\int_{-\infty}^{\infty}\left(\int_{-\infty}^{\infty}G(x,y)dx\right)dy=\int_{-\infty}^{\infty}\left(\int_{-\infty}^{\infty}G(x,y)dy\right)dx
\]
 which is the same as

\[
\int_{-\infty}^{\infty}\phi(y)\left(\int_{-\infty}^{\infty}f(x)e^{-iyx}dx\right)dy=\int_{-\infty}^{\infty}f(x)\left(\int_{-\infty}^{\infty}\phi(y)e^{-iyx}dy\right)dx.
\]
 If we divide both sides by $2\pi$ this is exactly (\ref{aswd}).
$\blacksquare$

It is the formula (\ref{aswd}) which will tell us how to define the
Fourier transform of elements of $\mathcal{S}^{\prime}$: If $\phi$
happens to be a function in $\mathcal{S}$ then this formula can be
rewritten as: 
\begin{equation}
L_{\widehat{f}}(\phi)=L_{f}(\widehat{\phi}).\label{jc}
\end{equation}
 It seems reasonable to think of $L_{\widehat{f}}$ as the Fourier
transform of $L_{f}$. So, to go one step further, if $L$ is some
other element of $\mathcal{S}^{\prime}$ which is not generated by
a function in $G(\Bbb{R)}$, we can just copy (\ref{jc}) and define
the Fourier transform of $L$ to be a new map from $\mathcal{S}$
to $\Bbb{C}$, which we will denote by $\widehat{L}$ and which acts
according to the formula 
\begin{equation}
\widehat{L}(\phi)=L(\widehat{\phi})\text{ for all }\phi\in\mathcal{S}.\label{defft}
\end{equation}

Here of course we need to know that $\widehat{\phi}\in\mathcal{S}$,
but we have already shown that, as property \ref{sfour}. We now want
to check that $\widehat{L}\in\mathcal{S}^{\prime}$. First, since
the Fourier transform of functions is linear, and $L$ is linear we
easily deduce that $\widehat{L}$ is also linear. Here is the proof:
\[
\widehat{L}(\alpha\phi+\beta\psi)=L\left(\widehat{\alpha\phi+\beta\psi}\right)=L\left(\alpha\widehat{\phi}+\beta\widehat{\psi}\right)=\alpha L(\widehat{\phi})+\beta L(\widehat{\psi})=\alpha\widehat{L}(\phi)+\beta\widehat{L}(\psi).
\]
 So now it remains to show that $\widehat{L}$ satisfies the second
``continuity'' property of Definition \ref{td}. This is an immediate
consequence of the following result:

\begin{lemma} \label{ftic}Let $\left\{ \phi_{k}\right\} _{k\in\Bbb{N}}$
be a sequence of functions in $\mathcal{S}$ which satisfies

\noindent ${\displaystyle \lim_{k\rightarrow\infty}C_{m,n}(\phi_{k})=0}$
for all non negative integers $m$ and $n$. Then ${\displaystyle \lim_{k\rightarrow\infty}C_{m,n}(\widehat{\phi_{k}})=0}$
for all non negative integers $m$ and $n$. \end{lemma}

The proof of Lemma \ref{ftic}, which will be deferred to an appendix
(Subsection \ref{pftic}), uses the same ideas as were used to prove
property \ref{sfour} of $\mathcal{S}$, (in particular (\ref{spz}),
(\ref{doft}) and (\ref{ftod})) with a more careful writing down
of estimates.

Now let us calculate the Fourier transform of some particular distributions.
For example, if $L=\delta_{0}$, then $\widehat{\delta_{0}}$ must
satisfy 
\[
\widehat{\delta_{0}}(\phi)=\delta_{0}(\widehat{\phi})=\widehat{\phi}(0)=\frac{1}{2\pi}\int_{-\infty}^{\infty}e^{-i0x}\phi(x)dx=\int_{-\infty}^{\infty}\frac{1}{2\pi}\phi(x)dx.
\]
 for all $\phi\in\mathcal{S}$. So we see that $\widehat{\delta_{0}}(\phi)=L_{g}(\phi)$
for all $\phi\in\mathcal{S}$, where $g$ is the constant function
$g(x)=\frac{1}{2\pi}$. This is the precise version of what we guessed
to be true (cf.\ (\ref{gad})) at the beginning of these notes.

More generally, suppose that $L=D^{n}\delta_{a}$ for some integer
$n\ge0$ and some $a\in\Bbb{R}$. Then, for every $\phi\in\mathcal{S}$,
we see that $\widehat{L}=\widehat{D^{n}\delta_{a}}$ must satisfy
\[
\widehat{D^{n}\delta_{a}}(\phi)=D^{n}\delta_{a}(\widehat{\phi})=(-1)^{n}\frac{d^{n}}{dx^{n}}\widehat{\phi}(a).
\]
 By (\ref{doft}) this equals $(-1)^{n}\cdot\frac{1}{2\pi}\int_{-\infty}^{\infty}e^{-iax}(-ix)^{n}\phi(x)dx=\int_{-\infty}^{\infty}\frac{1}{2\pi}e^{-iax}(ix)^{n}\phi(x)dx$.
This means that $\widehat{D^{n}\delta_{a}}$ is a function. More precisely,
it is the distribution $L_{f}$ corresponding to the function $f(x)=\frac{1}{2\pi}e^{-iax}(ix)^{n}$.

Next we shall calculate the Fourier transform of $x^{n}$, i.e.\ of
the distribution $L_{g}$ corresponding to the function $g(x)=x^{n}$
for some integer $n\ge0$. Perhaps you can already guess the answer
from the previous example, if you suppose that there is a connection
between Fourier transforms and inverse Fourier transforms of tempered
distributions.

For every $\phi\in\mathcal{S}$, we have, using (\ref{defft}), (\ref{ftod})
and then (\ref{ift}), that 
\begin{eqnarray*}
\widehat{L_{g}}(\phi) & = & L_{g}(\widehat{\phi})=\int_{-\infty}^{\infty}g(x)\widehat{\phi}(x)dx=\int_{-\infty}^{\infty}x^{n}\widehat{\phi}(x)dx=(-i)^{n}\int_{-\infty}^{\infty}(ix)^{n}\widehat{\phi}(x)dx\\
 & = & (-i)^{n}\int_{-\infty}^{\infty}\widehat{\phi^{(n)}}(x)dx=(-i)^{n}\int_{-\infty}^{\infty}\widehat{e^{i0x}\phi^{(n)}(x)}dx=(-i)^{n}\phi^{(n)}(0)\\
 & = & i^{n}D^{n}\delta_{0}(\phi).
\end{eqnarray*}

Since this is true for all $\phi\in\mathcal{S}$, we have just shown
that $\widehat{L_{x^{n}}}=i^{n}D^{n}\delta_{0}$.

Finally let us extend this last example by linearity to show that
the distributional Fourier transform of any polynomial $p(x)=\sum_{k=0}^{n}a_{k}x^{k}$
is given by $\sum_{k=0}^{n}a_{k}i^{k}D^{k}\delta_{0}$. An abbreviated
way of writing this might be {\large ``}$\widehat{p}=p(iD)\delta_{0}${\large ''}.
I enclosed this formula in quotation marks because both sides have
to be interpreted carefully.

\section{\label{conc}Some concluding remarks.}

\smallskip{}
These few pages, as we already said, can only give a very quick introduction
to tempered distributions. We shall conclude these notes by saying
something about notation, and then considering and offering some quick
and partial answers to two ``natural'' questions, an abstract one
and then a slightly more ``practical'' one.

\textit{0. About notation.}

In many books you will see that people write integrals with delta
functions as if they were ordinary integrals, i.e., instead of the
precise formula (in fact definition) that $\delta_{0}(\phi)=\phi(0)$
people like to write $\int_{-\infty}^{\infty}\delta_{0}(x)\phi(x)dx=\phi(0)$.
More generally, if $L$ is some general distribution in $\mathcal{S}^{\prime}$
then sometimes people like to pretend that it is somehow like an ordinary
function, and so, instead of the notation $L(\phi)$ which we have
used here, they use the ``integral'' $\int_{-\infty}^{\infty}L(x)\phi(x)dx$
to denote the value of $L$ when it acts on the function $\phi$.
This is fine and sometimes even useful, provided you remember that
it is \textit{only} notation, that it is \textit{not} a real integral,
and that in general, for individual values of $x$ the symbol $L(x)$
may be \textit{completely meaningless}. If you do not remember these
things then you can easily get to all sorts of impossible and illogical
conclusions, as we saw for example at the beginning of these notes.

Here is one example of how this ``integral'' notation can be used.
Suppose $L$ is a distribution and $\phi$ is a test function and
$c$ is a real constant. What is the exact meaning of the ``integral''
$\int_{-\infty}^{\infty}L(x+c)\phi(x)dx$ ? In other words, what is
the distribution $L(x+c)$? Well, we could guess that we would like
to have $\int_{-\infty}^{\infty}L(x+c)\phi(x)dx=$ $\int_{-\infty}^{\infty}L(x)\phi(x-c)dx$,
even though both of these things are not really integrals. But now
we see that a logical interpretation for both of them is to say that,
for each test function $\phi$, they equal $L(\psi)$ where $\psi$
is the test function defined by $\psi(x)=\phi(x-c)$ for all $x\in\Bbb{R}$.

\smallskip{}

\textit{1. How complicated and nasty can tempered distributions be?}

Well they cannot be too terrible. It turns out that every tempered
distribution can be constructed from a finite collection of continuous
functions in a finite number of steps involving distributional differentiation.
So we could say that tempered distributions are \textit{relatively}
simple objects. Here is a precise version of this result, the so called
``structure theorem''. (See pages 239--240 of \cite{Sc}for a slightly
different but equivalent formulation. For a more general version 
which applies to the case of functions
and distributions of $n$ variables, see, for example, Theorem 25.4
on pages 272--273 of \cite{T}.)

\textit{Given any tempered distribution }$L$\textit{, there exists
a finite collection of continuous functions }$f_{m}:\Bbb{R\rightarrow\Bbb{R}}$\textit{,
}$m=1,2,...,n$\textit{\ which each satisfy }$\left|f_{m}(x)\right|\le C\left(1+|x|\right)^{N}$\textit{\ for
some constant }$C$\textit{\ and some integer }$N$\textit{\ and
all }$x\in\Bbb{R}$\textit{\
and such that }$L$\textit{\ is given by } 
\[
L=\sum_{m=1}^{n}D^{m}L_{f_{m}}.
\]
 This formula is of course the same as the condition 
\[
L(\phi)=\sum_{m=1}^{n}(-1)^{m}\int_{-\infty}^{\infty}f_{m}(x)\frac{d^{m}}{dx^{m}}\phi(x)dx\text{ for all }\phi\in\mathcal{S}.
\]

\textit{2. What sort of things can distributions be used for?}

\smallskip{}
I answer this as a pure mathematician. But I hope that, with the help
of colleagues working in other fields, future versions of these notes
will also mention other applications, such as those encountered in
the courses which electrical engineering students take after they
complete our basic course about Fourier transforms.

One important use of distributions is in the study of differential
equations, including partial differential equations. In these notes
we only considered functions of one variable. If we wish to consider
partial differential equations we have to deal with functions of several
variables, and their analogous generalized functions. There is a natural
extension to $\Bbb{R}^{n}$ of the Fourier transform, the class $\mathcal{S}$,
and so also the class $\mathcal{S}^{\prime}$. Sometimes when it is
not at all clear that a differential equation has a solution, i.e.\ a
function satisfying the equation, it turns out that there are ways
to show that there is a distribution which satisfies the same equation.
Once it is known that there is a distributional solution, this can
sometimes be the first step towards showing that there is also a solution
in the original sense of the word. See \cite{Ru}, especially Chapter
8, for some examples of applications of distributions to partial differential
equations.

\section{\label{appendix}Appendices.}

\subsection{\label{csp}Functions in $\mathcal{S}$ which vanish outside a given
interval.}

\smallskip{}
The first and main step for constructing such functions is to consider
the function $v:\Bbb{R\rightarrow\Bbb{R}}$ given by 
\[
v(x)=\left\{ \begin{array}{l}
e^{-1/x^{2}}\text{ if }x>0\\
0\text{ if }\le0
\end{array}\right..
\]
 We will show that $v$ is infinitely differentiable on $\Bbb{R}$,
i.e.\ the $n^{th}$ derivative $v^{(n)}(x)$ of $v$ exists for each
$n\in\Bbb{N}$ and for all $x\in\Bbb{R}$. For all $x>0$ we have
$v^{\prime}(x)=\frac{2}{x^{3}}e^{-1/x^{2}}$ and $v^{\prime\prime}(x)=-\frac{6}{x^{4}}e^{-\frac{1}{x^{2}}}+\frac{4}{x^{6}}e^{-\frac{1}{x^{2}}}$.
Both of these derivatives are finite sums of functions of the form
$\frac{C}{x^{m}}e^{-1/x^{2}}$ where $C$ and $m\in\Bbb{N}$ are constants.
We do not need an explicit formula for $v^{(n)}$, but, by induction,
we can continue and show that, for all $x>0$ and each integer $n\ge0$,
$v^{(n)}(x)$ is a finite sum of functions of the form $\frac{C}{x^{m}}e^{-1/x^{2}}$.
Since $\lim_{x\rightarrow0+}\frac{1}{x^{m}}e^{-1/x^{2}}=\lim_{t\rightarrow\infty}t^{m}e^{-t^{2}}=0$
for each $m\ge0$, we see that $\lim_{x\rightarrow0+}v^{(n)}(x)=0$
and also that 
\begin{equation}
\lim_{x\rightarrow0+}\frac{1}{x}v^{(n)}(x)=0.\label{lfr}
\end{equation}
 In particular, for $n=1$, this second limit shows that $\lim_{h\rightarrow0+}\frac{v(h)-v(0)}{h}=0$.
Obviously also $\lim_{h\rightarrow0-}\frac{v(h)-v(0)}{h}=0$. So $v^{\prime}(0)$
exists and equals $0$. Obviously $v^{\prime}(x)$ also exists for
all $x\ne0$.

Now let us use induction. If we know that $v^{(n)}(x)$ exists for
all $x\in\Bbb{R}$ and also that $v^{(n)}(0)=0$, then (\ref{lfr})
tells us that $\lim_{h\rightarrow0+}\frac{v^{(n)}(h)-v^{(n)}(0)}{h}=0$
and, much as before, we deduce that $v^{(n+1)}(0)=0$ and $v^{(n+1)}(x)$
exists for all $x\in\Bbb{R}$.

Now, given any bounded interval $[a,b]$, let $\phi(x)=v(x-a)v(b-x)$
for all $x\in\Bbb{R}$. Since $v(x-a)$ and $v(b-x)$ are both infinitely
differentiable, so is $\phi(x)$ and obviously $\phi(x)=0$ for all
$x\le a$ and all $x\ge b$, and $\phi(x)>0$ for all $x\in(a,b)$.
It is obvious that $\phi\in\mathcal{S}$.

\subsection{\protect\smallskip{}
\label{recover}Recovering $f$ from $L_{f}(\phi)$.}

\begin{theorem} \smallskip{}
\label{rec}Let $\phi$ be a fixed even function in $\mathcal{S}$
such that $\int_{-\infty}^{\infty}\phi(x)dx=1$. For each $k\in\Bbb{N}$
and $c\in\Bbb{R}$ define $\phi_{k,c}$ by $\phi_{k,c}(x)=k\phi(k(x-c))$.
Then, for each $f\in SG(\Bbb{R)}$, 
\begin{equation}
f(c)=\lim_{k\rightarrow\infty}\int_{-\infty}^{\infty}f(x)\phi_{k,c}(x)dx.\label{apid}
\end{equation}
 \end{theorem}

In other words, we can determine the value of $f(c)$ at each point
$c\in\Bbb{R}$ if we know the values of $L_{f}(\phi_{k,c})$ for each
$k\in\Bbb{N}$. There are of course many ways to choose $\phi$. For
example we can take $\phi(x)=\frac{e^{-x^{2}}}{\sqrt{\pi}}$ or we
can choose $\phi(x)=cv(x+a)v(a-x)$ for some $a>0$, i.e.\ one of
the functions constructed in Subsection \ref{csp} multiplied by a
suitable constant $c$.

\textit{Proof. }Some parts of the proof of this theorem may perhaps
remind you of the proof in \cite{PZh} and \cite{PZe} of the inverse
Fourier theorem. The proof can be made rather simpler in the case
where the function $\phi$ vanishes outside some interval, but we
shall treat the general case. We first use the change of variables
$t=k(x-c)$ to obtain that

\begin{eqnarray*}
\int_{-\infty}^{\infty}f(x)k\phi(k(x-c))dx & = & \int_{-\infty}^{\infty}f\left(c+\frac{t}{k}\right)\phi(t)dt\\
 & = & \int_{-\infty}^{0}f\left(c+\frac{t}{k}\right)\phi(t)dt+\int_{0}^{\infty}f\left(c+\frac{t}{k}\right)\phi(t)dt.
\end{eqnarray*}
 We shall show that 
\begin{equation}
\lim_{k\rightarrow\infty}\int_{0}^{\infty}f\left(c+\frac{t}{k}\right)\phi(t)dt=\frac{1}{2}f(c+)\text{ and }\lim_{k\rightarrow\infty}\int_{-\infty}^{0}f\left(c+\frac{t}{k}\right)\phi(t)dt=\frac{1}{2}f(c-)\text{.}\label{dbb}
\end{equation}
 Since $f$ satisfies condition (\ref{avpm}) this will give (\ref{apid}).
We will only give the proof of the first part of (\ref{dbb}) since
the proof of the second part is almost exactly the same. Since $\phi$
is even and $\int_{-\infty}^{\infty}\phi(x)dx=1$, we have $\int_{0}^{\infty}\phi(t)dt=\frac{1}{2}$
and so the first part of (\ref{dbb}) is equivalent to the formula
\begin{equation}
\lim_{k\rightarrow\infty}\int_{0}^{\infty}\left(f\left(c+\frac{t}{k}\right)-f\left(c+\right)\right)\phi(t)dt=0.\label{pgmw}
\end{equation}

Given $\epsilon>0$, there exists $\delta=\delta(\epsilon)>0$ such
that $|f(c+s)-f(c+)|\le\epsilon$ for all $s\in(0,\delta)$. Bearing
this in mind we write 
\[
\left|\int_{0}^{\infty}\left(f\left(c+\frac{t}{k}\right)-f\left(c+\right)\right)\phi(t)dt\right|\le\int_{0}^{\infty}\left|\left(f\left(c+\frac{t}{k}\right)-f\left(c+\right)\right)\phi(t)\right|dt
\]
 
\[
\le\int_{0}^{\delta k}\left|\left(f\left(c+\frac{t}{k}\right)-f\left(c+\right)\right)\phi(t)\right|dt+\int_{\delta k}^{\infty}\left|\left(f\left(c+\frac{t}{k}\right)-f\left(c+\right)\right)\phi(t)\right|dt.
\]
 Then the first of the two integrals on the right is dominated by
\begin{equation}
\epsilon\int_{0}^{\delta k}|\phi(t)|dt\le\epsilon\int_{0}^{\infty}|\phi(t)|dt\label{ftt}
\end{equation}
 and the second integral is dominated by 
\begin{equation}
\int_{\delta k}^{\infty}\left|f\left(c+\right)\phi(t)\right|dt+\int_{\delta k}^{\infty}\left|f\left(c+\frac{t}{k}\right)\phi(t)\right|dt.\label{ebb}
\end{equation}

Since $\phi$ is absolutely integrable, the right side in (\ref{ftt})
is finite, and the first integral in (\ref{ebb}) equals $\left|f\left(c+\right)\right|\int_{\delta k}^{\infty}\left|\phi(t)\right|dt$
and tends to $0$ as $k$ tends to $\infty$. To estimate the second
integral we finally have to use the fact that $f$ satisfies the estimates
(\ref{vqp}) for some positive constants $C$ and $N$. We also make
two more changes of variable, first $y=t/k$, then later $z=c+y$.
This gives

\begin{eqnarray*}
\int_{\delta k}^{\infty}\left|f\left(c+\frac{t}{k}\right)\phi(t)\right|dt & = & k\int_{\delta}^{\infty}\left|f\left(c+y\right)\phi(ky)\right|dy\\
 & = & k\sum_{n=1}^{\infty}\int_{\delta n}^{\delta(n+1)}\left|f\left(c+y\right)\phi(ky)\right|dy\\
 & \le & k\sum_{n=1}^{\infty}\frac{C_{N+2,0}(\phi)}{(k\delta n)^{N+2}}\int_{\delta n}^{\delta(n+1)}\left|f\left(c+y\right)\right|dy\\
 & = & k\sum_{n=1}^{\infty}\frac{C_{N+2,0}(\phi)}{(k\delta n)^{N+2}}\int_{c+\delta n}^{c+\delta(n+1)}\left|f\left(z\right)\right|dz\\
 & \le & k\sum_{n=1}^{\infty}\frac{C_{N+2,0}(\phi)}{(k\delta n)^{N+2}}\int_{-|c|-\delta(n+1)}^{|c|+\delta(n+1)}\left|f\left(z\right)\right|dz\\
 & \le & k\sum_{n=1}^{\infty}\frac{C_{N+2,0}(\phi)}{|k\delta n|^{N+2}}C(1+|c|+\delta(n+1))^{N}\\
 & = & k^{-N-1}C_{N+2,0}(\phi)C\sum_{n=1}^{\infty}\frac{1}{(\delta n)^{N+2}}(1+|c|+\delta(n+1))^{N}.
\end{eqnarray*}
 Using easy estimates similar to those in (\ref{cls}) we see that
the general term in the series in this last expression is dominated
by a constant multiple of $\frac{1}{n^{2}}$. So the series converges
and $\int_{\delta k}^{\infty}\left|f\left(c+\frac{t}{k}\right)\phi(t)\right|dt\le k^{-N-1}M$
where the constant $M$ does not depend on $k$. Combining all these
estimates we obtain that

\[
\limsup_{k\rightarrow\Bbb{\infty}}\left|\int_{0}^{\infty}\left(f\left(c+\frac{t}{k}\right)-f\left(c+\right)\right)\phi(t)dt\right|\le\epsilon\int_{0}^{\infty}|\phi(t)|dt.
\]
 Then, since $\epsilon$ can be chosen arbitrarily small, this implies
that 
\[
\limsup_{k\rightarrow\Bbb{\infty}}\left|\int_{0}^{\infty}\left(f\left(c+\frac{t}{k}\right)-f\left(c+\right)\right)\phi(t)dt\right|=0
\]
 which gives (\ref{pgmw}) which, together with an analogous result
for

\noindent $\int_{-\infty}^{0}\left(f\left(c+\frac{t}{k}\right)-f\left(c-\right)\right)\phi(t)dt$
completes the proof of the theorem. $\blacksquare$

\subsection{\protect\smallskip{}
\label{delta}The distribution $\delta_{0}$ is not given by a function
in $SG(\Bbb{R)}$.}

This follows easily from Theorem \ref{rec}. Suppose, on the contrary,
that there exists a function $f\in SG(\Bbb{R)}$ such that $\delta_{0}=L_{f}$.
This means that $\psi(0)=\delta_{0}(\psi)=L_{f}(\psi)=\int_{-\infty}^{\infty}f(x)\psi(x)dx$
for all $\psi\in\mathcal{S}$. In particular, if we apply this to
the functions $\psi(x)=\phi_{k,c}(x)=k\phi(k(x-c))$ introduced in
Theorem \ref{rec} we obtain that $f(c)=\lim_{k\rightarrow\infty}\int_{-\infty}^{\infty}f(x)\phi_{k,c}(x)dx=\lim_{k\rightarrow\infty}\phi_{k,c}(0)=\lim_{k\rightarrow\infty}k\phi(-kc)$.
For each $c\ne0$, since $\left|\phi(-kc)\right|\le C_{1,0}(\phi)/|kc|$,
this last limit is $0$. So $f(c)=0$ for all $c\ne0$. For $c=0$
we may also obtain $f(0)=0$ if we choose a function $\phi$ satisfying
$\phi(0)=0$, or otherwise, if $\phi(0)$ is real, we will obtain
that $f(0)$ is $\infty$ or $-\infty$ depending on the sign of $\phi(0)$.
These conclusions contain several contradictions: First of all the
value of $f(0)$ as given by the formula (\ref{apid}) with $c=0$
should not depend on any particular choice of the function 
$\phi\mathbf{\in}\mathcal{S}$
which satisfies $\int_{-\infty}^{\infty}\phi(x)dx=1$. It is supposed
to be the same for all such $\phi$. Then, $f(0)$ should be finite
if $f\in SG(\Bbb{R)}$. Finally, if we are talking about integrals
in the usual sense of the word, changing the value of $f$ at $0$
should not influence the value of $\int_{-\infty}^{\infty}f(x)\phi(x)dx$
and this integral should thus equal $0$ for all choices of $\phi\in S$,
i.e.\ it will not in general equal $\phi(0)$. These contradictions
show that $\delta_{0}$ cannot equal $L_{f}$ for any $f\in SG(\Bbb{R)}$.

\subsection{\label{subseq}Integration by parts on an infinite interval.}

\smallskip{}
Let us complete the proof of (\ref{ibpp}) for all differentiable
$f\in SG(\Bbb{R)}$ and $\phi\in\mathcal{S}$ such that $f^{\prime}$
is also in $SG(\Bbb{R)}$. In view of property \ref{abint}, the functions
$f^{\prime}\phi$ and $f\phi^{\prime}$ are both absolutely integrable
on $\Bbb{R}$ and therefore also on $[0,\infty)$. So the 
limits $\lim_{r\rightarrow+\infty}\int_{0}^{r}f^{\prime}(x)\phi(x)dx$
and $\lim_{r\rightarrow+\infty}\int_{0}^{r}f(x)\phi^{\prime}(x)dx$
both exist. By (\ref{pib}), 
$f(r)\phi(r)=f(0)\phi(0)+\int_{0}^{r}f^{\prime}(x)\phi(x)dx
+\int_{0}^{r}f(x)\phi^{\prime}(x)dx$.
So, letting $r$ tend to $\infty$, we see 
that the limit $\lim_{r\rightarrow+\infty}f(r)\phi(r)$
exists. Now suppose that $\lim_{r\rightarrow+\infty}f(r)\phi(r)=c>0$.
Then for some sufficiently large $r$ we will have 
$f(x)\phi(x)\ge c/2$ for all $x\in[r,\infty)$. But this is impossible
since, again by property \ref{abint}, the function $f\phi$ is also
absolutely integrable. Similarly we get a contradiction if $c<0$.
It follows that $\lim_{r\rightarrow+\infty}f(r)\phi(r)=0$. Similarly
we can show that $\lim_{r\rightarrow-\infty}f(r)\phi(r)=0$ and so
(\ref{ibpp}) follows from (\ref{pib}).

(We don't have to worry about this here, but in fact, it can be shown
that $f^{\prime}$ has to be continuous on $\Bbb{R}$. This is because
of a certain property of functions which are derivatives of other
functions. Suppose that $g$ is a function which is the derivative
of some other function. There are examples which show that $g$ does
not have to be continuous. But $g$ can only be discontinuous in certain
ways. In particular, $g$ cannot have any simple jump discontinuities,
i.e, there are no points $x_{0}$ where $g(x_{0}+)$ and $g(x_{0}-)$
both exist and are different from each other.)

\subsection{\protect\smallskip{}
\label{name}A reason for using the notation $\mathcal{S}^{\prime}$.
Dual spaces, and continuous linear functionals.}

The use of the notation $\mathcal{S}^{\prime}$ is consistent with
notation which is often used in mathematics, in particular in the
field called \textbf{\textit{functional analysis}}. Suppose that $V$
is a vector space (of functions, or of some other objects) and we
have defined what we mean by convergent sequences in $V$. In particular
this means that we have a definition of what it means for a given
sequence $\{v_{k}\}_{k\in\Bbb{N}}$ to converge to the zero vector
in $V$. (Sometimes this definition of convergence is made with the
help of some norm on $V$. Sometimes, as in the case of $V=\mathcal{S}$
we prefer a different kind of definition.) . Then it is a standard
procedure to define a new space, which is called the \textbf{\textit{dual
space}} of $V$, and which is denoted by $V^{\prime}$, or sometimes
by $V^{*}$. This space consists of all ``continuous'' linear maps
$L:V\rightarrow\Bbb{C}$. (Or if $V$ is a vector space over $\Bbb{R}$
we will consider linear maps $L:V\rightarrow\Bbb{R}$) Here ``continuous''
means that, for every sequence $\{v_{k}\}_{k\in\Bbb{N}}$ which converges
to the zero element of $V$ we must have $\lim_{k\rightarrow\infty}L(v_{k})=0$.
The maps $L$ are often called \textbf{\textit{continuous linear functionals
on}} $V$. In this framework we can see that the tempered distributions
are exactly the continuous linear functionals on $\mathcal{S}$ and
the space of all these distributions is the dual space of 
$\mathcal{S}$. So it is appropriate to use the notation $\mathcal{S}^{\prime}$
for this space.

\subsection{\protect\smallskip{}
\label{pftic}The proof of Lemma \ref{ftic}.}

For all $\phi\in\mathcal{S}$ and all $\omega\in\Bbb{R}$ we have
\begin{eqnarray}
\left|\widehat{\phi}(\omega)\right| & \le & \frac{1}{2\pi}\int_{-\infty}^{\infty}|\phi(x)|dx=\frac{1}{2\pi}\int_{-\infty}^{\infty}\frac{1}{1+x^{2}}|(1+x^{2})\phi(x)|dx\notag\\
 & \le & \frac{1}{2\pi}\int_{-\infty}^{\infty}\frac{1}{1+x^{2}}\cdot\left(C_{0,0}(\phi)+C_{2,0}(\phi)\right)dx=\frac{1}{2\pi^{2}}\cdot\left(C_{0,0}(\phi)+C_{2,0}(\phi)\right).\label{uonk}
\end{eqnarray}
 Given any sequence $\left\{ \psi_{k}\right\} _{k\in\Bbb{N}}$ of
elements in $\mathcal{S}$ and any non negative integers $m$ and
$n$ we first apply (\ref{doft}) to obtain that $\frac{d^{n}}{d\omega^{n}}\widehat{\psi_{k}}(\omega)$
is the Fourier transform of $\xi_{k}(x)=(-ix)^{n}\psi_{k}(x)$ which
is also a function in $\mathcal{S}$. Then, by (\ref{ftod}), $(i\omega)^{m}\widehat{\xi_{k}}(\omega)$
is the Fourier transform of $\xi_{k}^{(m)}(x)$ which of course is
also in $\mathcal{S}$. Now, if we choose $\phi=\xi_{k}^{(m)}$, then
$\widehat{\phi}(\omega)=(i\omega)^{m}\widehat{\xi_{k}}(\omega)=(i\omega)^{m}\frac{d^{n}}{d\omega^{n}}\widehat{\psi_{k}}(\omega)$.
If we substitute in (\ref{uonk}) this gives 
\[
\left|\omega^{m}\frac{d^{n}}{d\omega^{n}}\widehat{\psi_{k}}(\omega)\right|\le\frac{1}{2\pi^{2}}\cdot\left(C_{0,0}(\xi_{k}^{(m)})+C_{2,0}(\xi_{k}^{(m)})\right).
\]
 Taking the supremum (or maximum) in this inequality as $\omega$
ranges over $\Bbb{R}$ gives 
\[
C_{m,n}(\widehat{\psi_{k}})\le\frac{1}{2\pi^{2}}\cdot\left(C_{0,0}(\xi_{k}^{(m)})+C_{2,0}(\xi_{k}^{(m)})\right).
\]

From this inequality it is clear that to complete the proof of Lemma
\ref{ftic} we have to show that the condition 
\begin{equation}
\lim_{k\rightarrow\infty}C_{m,n}(\psi_{k})=0\text{ for all non negative integers }m\text{ and }n,\label{wwag}
\end{equation}
 implies that $\lim_{k\rightarrow\infty}C_{0,0}(\xi_{k}^{(m)})$ and
$\lim_{k\rightarrow\infty}C_{2,0}(\xi_{k}^{(m)})$ are both $0$ for
all non negative integers $m$ and $n$. By Leibniz' formula, we have
that 
\begin{equation}
\xi_{k}^{(m)}(x)=\frac{d^{m}}{dx^{m}}\left[(-ix)^{n}\psi_{k}(x)\right]=\sum_{j=0}^{m}\binom{m}{j}\frac{d^{m-j}}{dx^{m-j}}(-ix)^{n}\cdot\psi_{k}^{(m)}(x),\label{dfg}
\end{equation}
 where, as usual, $\binom{m}{j}=\frac{m!}{(m-j)!j!}$. Each term $\frac{d^{m-j}}{dx^{m-j}}(-ix)^{n}$
is $0$ if $m-j>n$. Otherwise it equals $(-i)^{n}\frac{n!}{(n-m+j)!}x^{n-m+j}$.
So (\ref{dfg}) implies that 
\[
\left|\xi_{k}^{(m)}(x)\right|\le\sum_{j=\max\{0,m-n\}}^{m}\binom{m}{j}\cdot\frac{n!}{(n-m+j)!}\cdot\left|x^{n-m+j}\psi_{k}^{(m)}(x)\right|
\]
 and so, for any integer $\gamma\ge0$, we have 
\begin{eqnarray*}
\left|x^{\gamma}\xi_{k}^{(m)}(x)\right| & \le & \sum_{j=\max\{0,m-n\}}^{m}\binom{m}{j}\cdot\frac{n!}{(n-m+j)!}\cdot\left|x^{n-m+j+\gamma}\psi_{k}^{(m)}(x)\right|\\
 & \le & \sum_{j=\max\{0,m-n\}}^{m}\binom{m}{j}\cdot\frac{n!}{(n-m+j)!}\cdot C_{n-m+j+\gamma}(\psi_{k})
\end{eqnarray*}
 and so 
\begin{equation}
0\le C_{\gamma,0}\left(\xi_{k}^{(m)}\right)\le\sum_{j=\max\{0,m-n\}}^{m}\binom{m}{j}\cdot\frac{n!}{(n-m+j)!}\cdot C_{n-m+j+\gamma}(\psi_{k}).\label{rmwfg}
\end{equation}
 By (\ref{wwag}) the right side of (\ref{rmwfg}) tends to $0$ as
$k\rightarrow\infty$ for each fixed $m$, $n$ and $\gamma$. So
$\lim_{k\rightarrow\infty}C_{\gamma,0}\left(\xi_{k}^{(m)}\right)=0$
for each $\gamma$, in particular for $\gamma=0$ and $\gamma=2$,
which, as explained above, is exactly what we need to complete the
proof of Lemma \ref{ftic}. $\blacksquare$

\section{\label{symbols}A list of some of the symbols used in these notes}

(The meanings of some of the terminology used in the definitions of
these symbols are recalled in Section \ref{reminder}.)

$\bullet$ $\Bbb{R}$, $\Bbb{C}$, $\Bbb{N}$ and $\Bbb{Z}$. As usual
these are, respectively, the sets of all real numbers, all complex
numbers, all positive integers, and all integers.

$\bullet$ $\chi_{E}$. Let $E$ be any subset of $\Bbb{R}$, (in
many cases $E$ will be an interval). Then the \textbf{\textit{characteristic
function}} of $E$, sometimes also called the \textbf{\textit{indicator
function}} of $E$, is the function $\chi_{E}$ defined by $\chi_{E}(x)=\left\{ \begin{array}{ll}
1 & \text{ if }x\in E\\
0 & \text{if }x\notin E\text{.}
\end{array}\right.$.

$\bullet$ $G(\Bbb{R)}$. As in \cite{PZh} and \cite{PZe} and Definition \ref{GandFT} above, 
$G(\Bbb{R)}$
denotes the set of all functions $f:\Bbb{R\rightarrow\Bbb{C}}$ which
are piecewise continuous on each bounded subinterval of $\Bbb{R}$
and which are absolutely integrable on $\Bbb{R}$.

$\bullet$ $\widehat{f}$. The Fourier transform $\widehat{f}$ of
a function $f$ is defined by slightly different formulae in different
books. Here, as in \cite{PZh} and \cite{PZe}, we use the formula
$\widehat{f}(\omega)=\frac{1}{2\pi}\int_{-\infty}^{\infty}e^{-i\omega x}f(x)dx$.
Later $\widehat{f}$ is also defined when $f$ is a tempered distribution.

$\bullet$ $PC_{bj}$. The class of functions which are ``piecewise
continuous with balanced jumps''. See the definition given immediately
after (\ref{avpm}).

$\bullet$ $SG(\Bbb{R)}$. The set of functions of ``slow growth''.
See Definition \ref{sg}.

$\bullet$ $\mathcal{S}$. The Schwartz class of very smooth and very
rapidly decaying functions. See Definition \ref{sc}.

$\bullet$ $\mathcal{S}^{\prime}$. The space of tempered distributions.
See Definition \ref{td}.

\section{\label{reminder}Some reminders about some relevant mathematical
notions.}

Some of these notions and definitions can be found in the books \cite{PZh}
and \cite{PZe}, and some are introduced in these notes.

$\bullet$ \textbf{Piecewise continuity.}

Let $E$ be a bounded interval. A function $f:E\rightarrow\Bbb{C}$
is said to be \textbf{\textit{piecewise continuous on}} $E$ if it
is continuous at every point of $E$, or, in the worst case, there
are only finitely many points of $E$ at which $f$ is not continuous
and has a simple jump discontinuity. More precisely, we require the
one--sided limits from the left and right, $\lim_{x\rightarrow c-0}f(x)$
and $\lim_{x\rightarrow c+0}f(x)$ to exist at every interior point
$c$ of $E$ and to equal each other for all but at most finitely
many interior points of $E$. If $E$ contains its left endpoint $a$,
then we also require the existence of $\lim_{x\rightarrow a+0}f(x)$.
Similarly if $E$ contains its right endpoint $b$, then we require
the limit $\lim_{x\rightarrow b-0}f(x)$ to exist. (However $f(a)$
and $f(b)$ do not have to equal these limits.) We should also mention
other notation used for one sided limits: $f(c-)$, or $f(c-0)$ or
$\lim_{x\nearrow c}f(x)$ for limits from the left, and $f(c+)$,
or $f(c+0)$ or $\lim_{x\searrow c}f(x)$ for limits from the right.

If $E$ is an unbounded interval, for example if $E=\Bbb{R}$, then
we have adopted the same convention as in some books 
(cf. \cite{PZe} p. 94 or \cite{PZh}
p. 104). I.e., we define piecewise continuity of a function $f$ on
$E$ to mean that $f$ is piecewise continuous on every bounded subinterval
of $E$. So $f$ can possibly have infinitely many jumps altogether,
but each bounded interval can contain at most finitely many of those
jumps.

$\bullet$ \textbf{Riemann integrals of complex functions on a bounded
interval} $[a,b]$\textbf{.}

If we already know what we mean by the Riemann integral 
\footnote{
Here I will not repeat the standard definitions and properties of Riemann
integrals of real functions on a bounded interval $[a,b]$. But, if you do want to
revise these things and/or to see what I claim is an easier different way for defining
and studying Riemann integrals of one or more variables, then I invite you to look at \cite{C}.}%
$\int_a^bf(x)dx$
of a \textit{real} valued function $f$
on a bounded interval $[a,b]$, then it is very easy to also define and work with 
$\int_a^bf(x)dx$ when $f$ takes \textit{complex} values.
The easiest way for us to do this here
is to write each complex valued function in the form $f=u+iv$
where $u:[a,b]\rightarrow\Bbb{R}$ and $v:[a,b]\rightarrow\Bbb{R}$
are the real valued functions which are the real and imaginary parts
of $f$. Then our definition
will be
that $f$ is Riemann integrable on $[a,b]$ if and only
if $u$ and $v$ are both Riemann integrable on $[a,b]$ and that
$\int_{a}^{b}f(x)dx=\int_{a}^{b}u(x)dx+i\int_{a}^{b}v(x)dx$. Straightforward
(but boring) calculations then show that the standard formulæ\ $\int_{a}^{b}\alpha f(x)
+\beta g(x)dx=\alpha\int_{a}^{b}f(x)dx+\beta\int_{a}^{b}g(x)dx$
and $\left|\int_{a}^{b}f(x)dx\right|\le\int_{a}^{b}|f(x)|\ dx$ hold
for all integrable $f$ and $g$ and all constants $\alpha$ and $\beta$
also when some or all of them are complex valued.\footnote{ 
Here is a different but equivalent approach to all this: First we simply and
naturally extend
Riemann's original definition of his integral (which is \textit{not} the more frequently
used equivalent definition via upper and lower sums) to the case of complex valued
functions. Then the above statement which connects the integrability and integral
of a complex valued function with the integrability and integrals of its real and
imaginary parts becomes a theorem instead of
a definition.}%

\smallskip{}
$\bullet$ \textbf{Generalized Riemann integrals, and absolute integrability,
on unbounded intervals.}

Suppose that $f:[0,\infty)\rightarrow\Bbb{C}$ has the property that
its restriction to the interval $[0,R]$ is integrable for each constant
$R>0$. Then we can play two games:

(i) If the limit $\lim_{R\rightarrow\infty}\int_{0}^{R}f(x)dx$ and
is finite, then we say that $f$ has a \textbf{\textit{generalized
Riemann integral on}} $[0,\infty)$ which we denote and define by
$\int_{0}^{\infty}f(x)dx=\lim_{R\rightarrow\infty}\int_{0}^{R}f(x)dx$
.

\smallskip{}
(ii) If the limit $\lim_{R\rightarrow\infty}\int_{0}^{R}|f(x)|dx$
is finite%
\footnote{Since $\int_{0}^{R}|f(x)|dx$ is a non decreasing function of $R$
this limit always exists. But it may be infinite.%
} then we say that $f$ is \textbf{\textit{absolutely integrable on
}} $[0,\infty)$.

It is not difficult to prove that if $f$ is absolutely integrable,
then $f$ itself also has a generalized Riemann integral. (This is
similar to the proof that an absolutely convergent series is also
convergent.) But (analogously to series) the reverse implication is
untrue%
\footnote{Consider for example the function $f$ which equals $(-1)^{n}/n$
on the interval $[n-1,n)$ for all $n\in\Bbb{N}$.%
}.

On the interval $(-\infty,0]$ there are analogous definitions of
generalized Riemann integrals $\int_{-\infty}^{0}f(x)dx=\lim_{R\rightarrow+\infty}\int_{-R}^{0}f(x)dx$
and absolute integrability for functions $f:(-\infty,0]\rightarrow\Bbb{C}$
which are Riemann integrable on each \smallskip{}
interval $[-R,0]$ for each $R>0$.

Finally, if the function $f:\Bbb{R\rightarrow\Bbb{C}}$ is Riemann
integrable on $[-R,R]$ for each $R>0$, we can play three games:

(i) We say that $f$ has a \textbf{\textit{generalized Riemann integral
on}} $\Bbb{(}-\infty,\infty)$ if the integrals $\int_{0}^{\infty}f(x)dx$
and $\int_{-\infty}^{0}f(x)dx$ defined above both exist. Then we
define and denote this integral by $\int_{-\infty}^{\infty}f(x)dx=\int_{-\infty}^{0}f(x)dx+\int_{0}^{\infty}f(x)dx$.

(ii) We say that $f$ is \textbf{\textit{absolutely integrable on}}
$(-\infty,\infty)$ if $|f|$ has a generalized Riemann integral on
$(-\infty,\infty)$.

(iii) If the ``symmetric'' limit $\lim_{R\rightarrow+\infty}\int_{-R}^{R}f(x)dx$
exists and is finite then we call its value the \textbf{\textit{``Cauchy
principal value'' integral of }}$f$\textbf{\textit{\ on}} $(-\infty,\infty)$,
and use the notation $P.V.\int_{-\infty}^{\infty}f(x)dx=\lim_{R\rightarrow+\infty}\int_{-R}^{R}f(x)dx$
for this special kind of integral.

From what we have already said above, it is clear that the absolute
integrability of $f$ on $(-\infty,\infty)$, (which is equivalent
to absolute integrability on both $[0,\infty)$ and $(-\infty,0]$)
implies the existence of $\int_{-\infty}^{\infty}f(x)dx$, but that
this integral also exists for functions $f$ which are not absolutely
integrable. The existence of $\int_{-\infty}^{\infty}f(x)dx$ obviously
implies the existence of $P.V.\int_{-\infty}^{\infty}f(x)dx$ and
then these two integrals are equal. If $f(x)\ge0$ then $\int_{-\infty}^{\infty}f(x)dx$
exists if and only if $P.V.\int_{-\infty}^{\infty}f(x)dx$ exists.
But simple examples, e.g.\ $f(x)=x$ show that $P.V.\int_{-\infty}^{\infty}f(x)dx$
can sometimes exist when $\int_{-\infty}^{\infty}f(x)dx$ does not.

One reason why principal value integrals $P.V.\int_{-\infty}^{\infty}f(x)dx$
are important is because of the formula for inverting the Fourier
transform. We know that, if $f$ is sufficiently nice, then $\frac{f(x+)+f(x-)}{2}=\lim_{R\rightarrow\infty}\int_{-R}^{R}\widehat{f}(\omega)e^{i\omega x}d\omega$
and this last integral is of course just $P.V.\int_{-\infty}^{\infty}\widehat{f}(\omega)e^{i\omega x}d\omega$.
(The function $f:\Bbb{R\rightarrow\Bbb{R}}$ defined by $f(x)=\frac{x}{|x|}e^{-x}$
for all $x\ne0$ and $f(0)=0$ provides an example showing that in
general we cannot replace $\lim_{R\rightarrow\infty}\int_{-R}^{R}$
by $\int_{-\infty}^{\infty}$ in the Fourier inversion formula.)

$\bullet$ \textbf{Vector spaces or linear spaces.} A \textbf{\textit{vector
space over the complex field}} is a set $V$ of elements, which we
often call \textbf{\textit{vectors}}, on which we have defined two
operations. The first operation is often called ``addition of vectors''
and usually denoted by $+$ (even though this same sign $+$ can have
many other meanings in other contexts). The second operation is called
multiplication of vectors by ``scalars'', 
which in our case are
complex numbers%
\footnote{The set of complex numbers $\Bbb{C}$ is only one example of an algebraic
object called a \textbf{\textit{field,}}which we will not define
here. Other examples of fields are the set $\Bbb{R}$
of real numbers and the set $\Bbb{Q}$ of rational numbers.
Our definition here of a complex vector space
is a special case of the more general notion of a vector space over
a field, which makes sense for any choice of field.} and is usually denoted 
by simply writing the scalar in front of the
vector, i.e.\ $\lambda v$ for some scalar $\lambda$ and vector
$v$. These two operations have to satisfy a number of conditions,
associative commutative and distributive laws, existence of a zero
element, etc. We presume that you know these conditions well from a course
in linear algebra, so we will not list them here. You can see them
listed in many books, for example they are presented as nine items 
on pages 5 and 6 of \cite{PZe} and on 
page 11 of \cite{PZh}. We sometimes simply
say \textbf{\textit{vector space }}or \textbf{\textit{\ linear space
}}or\textbf{\textit{\ complex vector space }}instead of ``vector
space over the complex field''.

In these notes we encounter the following five examples of sets which
can each be shown to be (infinite dimensional) complex vector spaces.
They are
$G(\Bbb{R)}$, $PC_{bj}$, $SR(\Bbb{R)}$, $\mathcal{S}$ and 
$\mathcal{S}^{\prime}$. And of course $\Bbb{C}$ itself is also obviously 
a one dimensional complex vector space with respect to the usual 
operations of addition and multiplication of complex numbers.

The sets $G(\Bbb{R)}$, $PC_{bj}$, $SR(\Bbb{R)}$, $\mathcal{S}$
and $\mathcal{S}^{\prime}$ have some common features. Each one of
them is a set of some kind of complex valued functions $f$ defined
on some kind of set $\Gamma$ and the operations of vector addition
and multiplication by scalars are defined ``pointwise''. (For each
one of the sets $G(\Bbb{R)}$, $PC_{bj}$, $SR(\Bbb{R)}$ and $\mathcal{S}$
the set $\Gamma$ is simply $\Bbb{R}$. In the more exotic case of
$\mathcal{S}^{\prime}$ we have to take $\Gamma=\mathcal{S}$.) This
means that six of the required above mentioned nine properties 
(on pages 5 and 6 of \cite{PZe} or page 11 of \cite{PZh})
follow immediately from analogous properties of multiplication and
addition of complex numbers.

Let me explain this in a bit more detail: Suppose that $V$ is one
of the above five sets and that $f$ and $g$ are two vectors in $V$.
Then they are both functions on $\Gamma$, i.e.\ for each $\gamma\in\Gamma$
we have complex numbers $f(\gamma)$ and $g(\gamma)$ which are the
values of $f$ and $g$ respectively at the ``point'' $\gamma$.
We define the new vector $f+g$ to be the function defined by the
formula $(f+g)(\gamma)=f(\gamma)+g(\gamma)$ for all $\gamma\in\Gamma$.
Note that in this last formula the symbol $+$ has two different meanings.
On the left side it means the operation of addition in $V$, and on
the right side it means the operation of addition of complex numbers.
Similarly, for each $\lambda\in\Bbb{C}$, we define $\lambda f$ to
be the function defined by $(\lambda f)(\gamma)=\lambda f(\gamma)$
for each $\gamma\in\Gamma$. (Here again, note that the writing of
$\lambda$ next to the vector $f$ and the writing of $\lambda$ next
to the number $f(\gamma)$ have two different meanings.) Now we want
to check that the operations defined in this way satisfy those above mentioned nine
conditions (on pages 5 and 6 of \cite{PZe} or on page 11 of \cite{PZh}). First we should consider conditions
1 and 5 which state that, whenever $f$ and $g$ are both in $V$
and $\lambda\in\Bbb{C}$, then $f+g$ and $\lambda f$ both have to
be in $V$. To verify these we need to use some special properties,
depending on the particular choice of $V$. E.g., we know that sums
of continuous functions are continuous, and sums of integrable functions
are integrable, sums of linear maps are linear, etc.\
etc. We might also use the triangle inequality for complex numbers
to show that certain required inequalities hold. Then condition 3
requires the existence of a special so--called zero element $\mathbf{0}$
in $V$ with the property that $f+\mathbf{0}=f$ for each f $f\in V$.
In each case that we have to consider, our set $V$ contains the zero
function, i.e the function whose value is $0$ (the complex number)
for each $\gamma\in\Gamma$. This clearly has the required property.

The verifications of the remaining six conditions are, as I said above,
almost automatic. For example, condition 2 is the associative law,
that $(f+g)+h=f+(g+h)$ for all $f$, $g$ and $h$ in $V$. In our
case this is simply the corresponding associative law for addition
of complex numbers applied to the numbers $f(\gamma)$, $g(\gamma)$
and $h(\gamma)$ for each $\gamma\in\Gamma$.

$\bullet$ \textbf{Separating classes} (or \textbf{total sets)}. Let
$A$ and $B$ be classes of functions on $\Bbb{R}$ such that the
integral $\int_{-\infty}^{\infty}f(x)g(x)dx$ exists for each $f\in A$
and each $g\in B$. Suppose that whenever $f_{1}$ and $f_{2}$ are
functions in $A$ such that $\int_{-\infty}^{\infty}f_{1}(x)g(x)dx=\int_{-\infty}^{\infty}f_{2}(x)g(x)dx$
for all $g\in B$, it follows that $f_{1}=f_{2}$. Then we say that
$B$ is a \textbf{\textit{separating class}} for $A$. (Or we can
also say that $B$ is a \textbf{\textit{total set }}for $A$.)

\phantom{$\bullet$}

\section{\protect\smallskip{}
\label{bib}Some books and papers for further reading.}

\smallskip{}
Here are some books and papers, most of which I happened to notice in 
our mathematics library. I have surely missed some other very good 
items, and will be glad to add them to this list if you tell me about 
them.

\bigskip{}

\phantom{.}
\end{document}